
\magnification=\magstep1
\input amstex
\documentstyle{amsppt}

\NoBlackBoxes
\pagewidth{5.75truein}
\pageheight{8.5truein}


\rightheadtext{Quadratic Linear Algebras}
\leftheadtext{(I. Gelfand, V. Retakh, R.L. Wilson)}

\topmatter
\title Quadratic linear algebras associated with 
factorizations of noncommutative polynomials and 
noncommutative differential polynomials 
\endtitle
\author Israel Gelfand, Vladimir Retakh, and Robert Lee Wilson
\endauthor
\address
\newline
Department of Mathematics, Rutgers University, Piscataway,
NJ 08854-8019
\endaddress

\email
\newline
I.~G. : igelfand\@ math.rutgers.edu
\newline
V.~R. : vretakh\@ math.rutgers.edu
\newline
R.~W. : rwilson\@math.rutgers.edu
\endemail
\leftheadtext{I.M. Gelfand, V. Retakh, R.L. Wilson}
\endtopmatter
{\bf Abstract}.
We study certain quadratic and quadratic linear algebras 
related to factorizations of noncommutative polynomials 
and differential polynomials. 
Such algebras possess a natural derivation 
and give us a new understanding of the 
nature of noncommutative symmetric functions.   
\bigskip 
\centerline{\bf Introduction}
\medskip
Let $x_1,\dots , x_n$ be the roots of a {\it generic} polynomial
$P(x)=x^n+a_1x^{n-1}+\dots +a_n$ over a division algebra $R$.
There are two important classical problems: a) to express the
coefficients $a_1,\dots , a_n$ through the roots, and b)
to determine all factorizations of $P(x)$, or $P(t)$, where
$t$ is a formal variable commuting with elements of $R$.

The first problem was solved in \cite {GR3, GR5}.
For any ordering $\{i_1,\dots , i_n\}$ of
$\{1,\dots , n\}$ elements $x_{\emptyset , i_1}=x_{i_1}$,
$x_{\{i_1, i_2, \dots ,i_{k-1}\}, i_k}\in R$, 
$k=2,\dots , n$, were constructed such that
for every $m=1, \dots, n$, 
$$
(-1)^ma_m=\sum _{j_1>j_2>\dots j_m}y_{j_1}y_{j_2}
\dots y_{j_m}, \tag 0.1
$$
where $y_1=x_{i_1}$, $y_k=x_{\{i_1,\dots , i_{k-1}\}, i_k}$, 
$k=2,\dots , n$.

It is surprising that the left-hand side in formula (0.1)
does not depend on the ordering of $\{1,\dots , n\}$ but the
right-hand side {\it a priori} depends on the ordering.
The independence of the right-hand side in (0.1) of the ordering of
$\{1,\dots , n\}$ was a key point in the theory of noncommutative
symmetric functions developed in \cite {GR3, GR5}.

The element  $x_{\{i_1,\dots , i_{k-1}\}, i_k}$ has an interesting structure.
It is symmetric in $x_{i_1}, \dots , x_{i_{k-1}}$.
It is a rational function in
$x_{i_1}, \dots , x_{i_k}$ containing $k-1$ inversions in the
generic case. In other words, $x_{\{i_1,\dots , i_{k-1}\}, i_k}$
is a rational expression of height $k-1$. 

In fact these elements satisfy simple relations:
$$
x_{A\cup \{i\},j}+ x_{A, i}=x_{A\cup \{j\},i}+ x_{A, j}, \tag 0.2a
$$
$$
x_{A\cup \{i\},j}\cdot x_{A, i}=x_{A\cup \{j\},i}\cdot x_{A, j}, \tag 0.2b 
$$
for all $A\subseteq \{1,\dots , n\}$, $i,j\notin A$.

In this paper, to avoid inversions, we define an algebra $Q_n$
to have generators $z_{A, i}$, 
$A\subseteq \{1,\dots , n\}$, $i\notin A$, and relations corresponding to
(0.2) (with $x$ replaced by $z$).

The algebra $Q_n$ is also a universal algebra for all possible
factorizations of $P(t)$. Set $A_k=\{i_1, \dots , i_{k-1}\}$
for $k=2,\dots , n$.
The formulas (0.1) are equivalent to the decomposition
$$
P(t)=(t-x_{A_n, i_n})(t-x_{A_{n-1}, i_{n-1}})\dots (t-x_{i_1}).
$$

In the generic case the polynomial $P(t)$ has $n!$ decompositions
into linear factors. We study all such
decompositions together using the algebra $Q_n$.

In this paper we study the internal structure of algebra $Q_n$.
It has a natural derivation and a natural anti-involution. 
For each pair $A, B\subseteq \{1,\dots , n\}$, with
$A\cap B=\emptyset $ we define an element $z_{A, B}\in Q_n$.
When $B$ contains more than one element 
these elements are ``invisible" in the commutative case, i.e.,
under the natural commutative specialization of $Q_n$
their image is zero. We are going to construct a
``noncommutative logic" using these elements.

Also, the elements $z_{A, \emptyset}$, and similarly, the elements
$z_{\emptyset, A}$,
for all $A\subseteq \{1,\dots , n\}$, $A\neq \emptyset $, constitute a basis
for the subspace of $Q_n$ spanned by all generators. These
elements satisfy simple quadratic relations.

One of our main results, Theorem 1.3.3 (see also Theorem 1.3.8),
is the determination of a basis for $Q_n$. 

For each ordering $I=(i_1,\dots , i_n)$ of $\{1,\dots , n\}$ there is
a natural subalgebra, denoted $Q_{n, I }$, of $Q_n$
generated by $\{z_{\{i_1, \dots , i_k\}, 
\emptyset } \ |\ 1\leq k\leq n\}$. Using the basis theorem
we describe arbitrary intersections of the $Q_{n, I}$.
In particular, the intersection of all $Q_{n, I }$ is
the algebra generated by the coefficients of $P(t)$.
 
The algebras $Q_n$ seem to be very interesting. 
They may be viewed as a Galois type extension of the 
algebra of symmetric functions.
A similar theory related to factorizations of differential 
polynomials over a noncommutative algebra is also presented in 
this paper.  Such factorizations 
(Miura decompositions) appeared in \cite{EGR} and were heavily 
used in the study of noncommutative integrable systems. 

The theory of factorizations of noncommutative polynomials and 
differential polynomials presented here is based on 
the theory of quasideterminants \cite {GR1-GR5}.
We recall the definition of quasideterminants.

Let $X=(x_{ij})$ be an $m\times m$-matrix over 
a division algebra $R$. For any $1\le i,j\le m$, let 
$r_i(X)$, $c_j(X)$ be the $i$-th row and the $j$-th column of $X$.  
Let $X^{ij}$ be the submatrix of $X$ obtained by removing 
the $i$-th row and the $j$-th column from $X$. 
For a row vector $r$ let $r^{(j)}$ be $r$ without 
the $j$-th entry.  For a column vector $c$ let $c^{(i)}$ 
be $c$ without the $i$-th entry. 
Assume that $X^{ij}$ is invertible. Then the quasideterminant 
$|X|_{ij}\in R$ is defined by the formula
$$
|X|_{ij}=x_{ij}-r_i(X)^{(j)}(X^{ij})^{-1}c_j(X)^{(i)},
$$
where $x_{ij}$ is the $ij$-th entry of $X$.

We thank S. Gelfand, F. Knop, S. Sahi and S. Serconek for
helpful remarks.

The second author was partially supported by the
National Science Foundation.
\bigskip

\head 1. The quadratic algebras $Q_n$\endhead

\subhead {1.1. Definition}
\endsubhead 
In this section we introduce
and study a quadratic algebra $Q_n$ over a field $k$ 
generated by formal variables
$z_{A,i}$ where $A\subset \{1,\dots ,n\}$ is an unordered set and
$i\in \{1,\dots ,n\}$, $i\notin A$. (The set $A$ might be
empty.) We also define a derivation
$\partial $ of $Q_n$, i.e., an endomorphism of
$Q_n$ such that 
$\partial (ab)=(\partial a)b + a\partial b$.

Let $F_n$ denote the free associative algebra over a field
$k$ generated by elements
$z_{A,i}$, $A\subset \{1,\dots ,n\}$, $i\in \{1,\dots , n\}$,
$i\notin A$. Let $J_n$ denote the ideal of $F_n$ generated by
the elements
$$z_{A\cup i, j}+z_{A,i}-z_{A\cup j, i}-z_{A,j}, \tag 1.1a $$
$$z_{A\cup i, j}\cdot z_{A,i}-
z_{A\cup j, i}\cdot z_{A,j}. \tag 1.1b $$
(We write $A\cup i$ and $A\setminus j$ instead of
$A\cup \{i\}$ and $A\setminus \{j\}$.)

Set $Q_n=F_n/J_n$ and denote the coset of $z_{A,i}$ in
$Q_n$ by the same expression $z_{A,i}$.
Clearly $S_n$, the symmetric group on $\{1,\dots , n\}$,
acts on $Q_n$ by $\sigma (z_{A,i})=z_{\sigma (A), \sigma (i)}$.
Note that the free algebra $F_n$ has 
a derivation $\partial $ defined by 
$$ \partial (z_{A,i})=1 \tag 1.1c $$
and an antiautomorphism $\theta $ defined by
$$ \theta (z_{A,i})=
z_{\{1,\dots ,n\}\setminus A\setminus i,i} \tag 1.1d $$
for all $A$, $i\notin A$.

\proclaim {Proposition 1.1.1} (a) The derivation $\partial $ 
of $F_n$ preserves $J_n$ and so induces a
derivation, again denoted $\partial $, of
$Q_n$ satisfying
$ \partial (z_{A,i})=1$ for all $A$, $i\notin A$.

(b) The antiautomorphism $\theta $ of $F_n$ preserves $J_n$
and so induces an antiautomorphism, again denoted $\theta $,
of $Q_n$ satisfying
$ \theta (z_{A,i})=
z_{\{1,\dots ,n\}\setminus A\setminus i,i} $
for all $A$, $i\notin A$.
\endproclaim

The proof of (a) follows from the fact
that the map $\partial $ applied to the relation (1.1b)
gives the relation (1.1a). This may be stated as follows: $Q_n$
is the differential algebra on generators $\{z_{A, i}\}$
defined by relations (1.1b) and (1.1c).
\smallskip
\noindent {\bf Example}. By definition $Q_2$ is
generated by the elements
$z_{\emptyset, i}$ (denoted by $z_i$) for
$i=1,2$ and  $z_{\{i\},j}$
(denoted by $z_{i,j}$) for $i=1, j=2$ and
$i=2, j=1$. These elements satisfy the relations
$$z_{1,2}+z_1=z_{2,1}+z_2, $$
and
$$z_{1,2}\cdot z_1=z_{2,1}\cdot z_2.$$

Thus $Q_2$ is generated by any three element subset of
$\{z_1, z_2, z_{1,2}, z_{2,1}\}$.
A more interesting choice of generators is
$$
\Lambda =z_1+z_{1,2}, \ \ u=z_{1,2}-z_2, \ \
\xi = z_1-z_2.
$$
Note that $\Lambda $ and $u$ are symmetric and $\xi$ is
skew-symmetric and that
$$
[\Lambda , \xi ]+[u, \xi]_+=0,
$$
where $[\ ,\ ]_+$ is the anticommutator.
For the anti-isomorphism $\theta $ one has
$\theta (z_1)=z_{2,1}$, $\theta (\Lambda )=\Lambda$,
$\theta (u)=-u$, $\theta (\xi )=\xi$.
 
In Section 2 we construct a natural map of $Q_n$ into 
the free skew-field generated by $n$ elements
$z_1,\dots , z_n$. If $\bar z$ is the image of
$z\in Q_n$ under this map, then  for each
$A$ and $i\notin A$
$$
\bar z_{A\cup i, j}=
(\bar z_{A,j}-\bar z_{A, i})\bar z_{A,j}(\bar z_{A, j}-
\bar z_{A, i})^{-1};
$$
this depends on the fact that 
$(z_{A, j}-z_{A, i})\neq 0$ in $Q_n$
(see Section 2) and on the following result.

\proclaim {Proposition 1.1.2} Let $A\subset \{1,\dots ,n\}$,
$i,j\in \{1,\dots ,n\}$, $i,j\notin A$. Then
$$
z_{A\cup i, j}(z_{A, j}-z_{A, i})=
(z_{A,j}-z_{A, i})z_{A,j}.
$$
\endproclaim

We now define two important homomorphisms of $Q_n$ into
commutative algebras. As usual, we let
$k[v_1,\dots , v_n]$ denote the (commutative) polynomial
algebra in $v_1,\dots , v_n$. Let $I_n$ denote the ideal
in $k[v_1,\dots , v_n]$ generated by $\{v_i^2+v_i,
1\leq i\leq n\}$. Let $K_n$ denote $k[v_1,\dots , v_n]/I_n$
and let $w_i=v_i+I_n$. Then there is a homomorphism
$$\phi : Q_n \to K_n $$
defined by 
$$\phi (z_{A,i})=w_i\prod _{j\in A}(1+w_j)$$
for all $A\subseteq \{1,\dots , n\}$, $i\notin A$,
and a homomorphism
$$\psi : Q_n \to k[v_1,\dots , v_n]$$
defined by
$$ \psi (z_{A,i})=v_i$$
for all $A\subseteq \{1,\dots , n\}$, $i\notin A$.

Let $T_n$ denote the span of $\{z_{A,i}\ |\
A\subseteq \{1,\dots , n\}, i\notin A\}\subseteq Q_n$
and $\bar T_n=T_n + k\cdot 1$.

\proclaim{Lemma 1.1.3} $\phi |_{\bar T_n} : \bar T_n \to K_n$ 
is an isomorphism of vector spaces.
\endproclaim

{\it Proof}. $\phi |_{\bar T_n}$ is an epimorphism
and $\text {dim}\ \bar T_n =\text {dim}\ K_n = 2^n$. 
 
\medskip

\subhead {1.2. The elements $z_{A,B}$ in $Q_n$}
\endsubhead

We define now elements $z_{A,B}\in Q_n$ for any
$A,B\subseteq \{1,\dots ,n\}$. We will show
that the elements $z_{A, \emptyset }$ and $z_{\emptyset , B}$
form two bases for the linear envelope of all
generators $z_{A, i}$.

For $A,B\subseteq \{1,\dots ,n\}$ define  $z_{A, B}$ to be the
projection of 
$$
(\phi |_{\bar T_n})^{-1}(\prod _{i\in B}w_i)(\prod _{j\in A}(1+w_j))
$$
on $T_n$.

\proclaim{Proposition 1.2.1} 

i) $\ \ z_{\emptyset, \emptyset}=0$, 

ii) $\ \ z_{A, B}=0$ if $A\cap B\neq \emptyset $, 

iii) $\ \ z_{A\cup i,B}-z_{A,B\cup i}=z_{A,B},\ \ i\notin A,B$,

iv) $\ \ z_{A,B\cup i}=\sum _{D\subseteq B}
(-1)^{|B|-|D|}z_{A\cup D,i}, \ \ i\notin A,B$.
\endproclaim

\proclaim{Corollary 1.2.2} i) If $A=\{i_1,\dots , i_r\}$,
then
$$
z_{A, \emptyset }=\Lambda _1(A)=z_{i_1}+
z_{i_1, i_2}+\dots z_{i_1\dots i_{r-1}, i_r}.
$$
ii) The elements $z_{A, \emptyset}$ and $z_{\emptyset, B}$
are connected by a M\"obius transformation
$$ z_{\emptyset , B}=\sum _{D\subseteq B}(-1)^{|B|-|D|}
z_{D,\emptyset},$$
$$z_{A,\emptyset }=\sum _{C\subseteq A}z_{\emptyset, C}.$$
iii) If $\sigma (A)=A$ and $\sigma (B)=B$, then $\sigma (z_{A, B})=
z_{A, B}$.
\endproclaim      

\noindent {\bf Example}.
$$
z_{\emptyset, 1}=z_1, \ \ 
z_{\emptyset , 12}=z_{1,2}-z_2=z_{2,1}-z_1,
$$
$$
z_{1, \emptyset }=z_1, \ \ 
z_{12,\emptyset }=z_1+z_{1,2}=z_2+z_{2,1}.
$$

Denote $z_{A, \emptyset }$ by $r(A)$ and 
$z_{\emptyset , B}$ by $u(B)$. We may express $z_{A,B}$
via those elements.

\proclaim{Proposition 1.2.3} Let $i\notin A$. Then
$$
z_{A,i}=r(A\cup i)-r(A),
$$
$$
z_{A,i}=\sum _{i\in D\subseteq A\cup i} u(D).
$$
\endproclaim

This proposition follows from the next statement.

\proclaim {Proposition 1.2.4} Let $A\cap B=\emptyset $. Then
$$
z_{A,B}=z_{A\cup \{i_1,\dots , i_l\}, B}-\sum _{k=0}^{l-1} z_{A\cup
\{i_1,\dots , i_k\}, B\cup i_{k+1}}
$$
for any $\{i_1,\dots , i_l\}\subseteq \{1,\dots , n\}\setminus (A\cup B)$,
and
$$
z_{A, B}=\sum _{B\subseteq D\subseteq (A\cup B)}u(D).
$$
\endproclaim

\proclaim {Proposition 1.2.5} If $|B|\geq 2$, then
$\psi (z_{A,B})=0$. 
\endproclaim

This shows that the elements $z_{A, B}$ for $|B|\geq 2$
carry the ``noncommutative structure" of the algebra $Q_n$.
Let $\partial $, $\theta $ be the derivation and the anti-isomorphism
defined in Section 1.1.

\proclaim{Proposition 1.2.6} 
$$
\partial (z_{A,B})=0\ \ \ 
{\text {if}}\ \ |B|\geq 2,
$$

$$
\theta (z_{A,B})=(-1)^{|B|+1}
z_{\{1\dots n\}\setminus (A\cup B), B}.
$$
\endproclaim
 
\proclaim {Proposition 1.2.7} Each of the families 
$\{z_{A,\emptyset }\}$, $\{z_{\emptyset,A}\}$
for all nonempty $A\subset \{1,\dots ,n\}$ forms a basis for
$T_n$. Another basis in $T_n$ is given by the elements
$z_{A,\bar A}$, where $\bar A$ is $\{1,\dots , n\} \setminus A$.
\endproclaim

Note, that 
$\sigma (z_{A,\bar A})=z_{A,\bar A}$ if $\sigma (A)=A$. 
 
The proof of Proposition 1.2.7 follows from Proposition 1.2.8
and Lemma 1.1.3.

\proclaim {Proposition 1.2.8} The elements $z_{A,\bar A}$ satisfy the
following formulas:
$$z_{C , D}=\sum _{C\subseteq A, A\cap D=\emptyset} 
(-1)^{|\bar A|-|D|}z_{A,\bar A},$$
$$ z_{A,\bar A}=\sum _{A\subseteq C}(-1)^{n-|C|}z_{C,
\emptyset },$$
$$z_{A,\bar A}=
\sum _{\bar A\subseteq D} z_{\emptyset ,D}.
$$
\endproclaim

\medskip                       

\subhead {1.3. Multiplicative relations and linear bases}
\endsubhead

We describe below the {\it multiplicative relations} for the elements
$z_{A,\emptyset }=r(A)$, $z_{\emptyset, A}=u(A)$
and their corollaries. 
The multiplicative relations for other
bases in $T_n$ can be written in a similar manner.

The relations (1.1) imply:
\proclaim {Proposition 1.3.1} Let $i,j\notin A$, $i\neq j$. 
Then
$$
(r(A\cup i\cup j)-r(A\cup i))\cdot (r(A\cup i)-r(A\cup j))=
$$
\noindent {\rm (1.2)}
$$
=(r(A\cup i)-r(A\cup j))\cdot (r(A\cup j)-r(A)),
$$
$$
[z_{A,i}, z_{A,j}]=
\sum _{\{ij\}\subseteq D\subset A\cup \{ij\}} 
u(D)\cdot (z_{A,i}-z_{A,j}). \tag 1.3
$$
\endproclaim

The relations (1.2) have a simple matrix form. We will not
use this form in this paper and so the reader can skip the
next proposition. For each
$j, 1\leq j \leq n$, define matrices $R(j)$ and $S(j)$ with
rows and columns indexed by $B,C\subset \{1,\dots ,n\}$,
$j\in B, C$, as follows:
$$
R(j)_{B,C}=\delta _{B,C}(r(B)-r(B\setminus j)),
$$
$$
S(j)_{B,B\setminus i}=r(B\setminus j)-r(B\setminus i)
$$
if $i\neq j$, $S_{B,C}(j)=0$ otherwise.
Note that $R(j)$ is a diagonal matrix.

Define $R$ to be the block diagonal matrix with blocks
$R(1),\dots ,R(n)$, $S$ to be the block diagonal matrix
with blocks $S(1),\dots ,S(n)$.

\proclaim {Proposition 1.3.2} The matrices $R$ and $S$ commute
if and only if the relations (1.2) are fulfilled.
\endproclaim

We will construct now a basis for $Q_n$. This construction is
based on the standard ordering $1<2<\dots <n$ of 
$\{1,\dots , n\}$.

A {\it string} is a finite sequence $\Cal B=(B_1,\dots, B_l)$
of nonempty subsets of $\{1,\dots ,n\}$. 
Let $A\subseteq \{1,\dots ,n\}$, $|A|=u$. Write
$A=\{a_1,\dots , a_u\}$ where $a_1>a_2>\dots >a_u$.
For $1\leq j <u$ define $(A:j)$ to be the string
$(A, A\setminus \{a_1\}, \dots , A\setminus \{a_1,
\dots , a_{j-1}\}$. 
Write $r(A:j)=r(A)r(A \setminus \{a_1\})...r(A\setminus \{a_1,
\dots , a_{j-1}\})$.  Then we have:

\proclaim {Theorem 1.3.3} The set of all products
$r(A_1:j_1)\dots r(A_s:j_s)$ where \break $A_1,\dots , A_s
\subseteq \{1,\dots , n\}$, $j_i\leq |A_i|$ for all $i$,
and, for all $2\leq i\leq s$, we have either 
$|A_i|\neq |A_{i-1}|-j_{i-1}$ or $A_i\nsubseteq A_{i-1}$
is a basis for $Q_n$. 
\endproclaim

The proof follows from Theorem 1.3.8 below. To formulate
Theorem 1.3.8 we need some definitions and notations.

Let $\Cal B =(B_1, \dots , B_l)$. We call
$l=l(\Cal B)$ the {\it length} of $\Cal B$ and 
$|\Cal B|=\sum _{i=1}^l |B_i|$ the {\it degree} of $\Cal B$.

If $\Cal B=(B_1, \dots , B_l)$ let $r(\Cal B)\in Q_n$
denote the product $r(B_1)\dots r(B_l)$.
For any set $Z$ of strings we will denote $\{r(\Cal B)|
\Cal B \in Z\}$ by $r(Z)$.  Define an increasing filtration $$Q_{n,0} 
\subseteq Q_{n,1} \subseteq ... $$ by $$Q_{n,i} = 
{\text {span}} \{r(\Cal B)| \ |\Cal B| \le i\}.$$

For $0 \le i \le l(\Cal B)$ denote the {\it truncated} string $(B_{i+1},\dots,B_l)$ by
$$T_i(\Cal B) = (T_i(\Cal B)_1 ,\dots , T_i(\Cal B)_{l-i}).$$  Thus 
$T_i(\Cal B)_j = B_{i+j}.$

Suppose $\Cal B=(B_1,\dots, B_l)$ is a string.
We will define by induction a sequence of integers
$n(\Cal B)=(n_1,n_2,\dots , n_t)$,
$1=n_1<n_2<\dots <n_t=l+1$, as follows:

$n_1=1$,
$n_{k+1}=\min (\{j>n_k|\ B_j\nsubseteq B_{n_k}\  {\text or}\ 
|B_j|\neq |B_{n_k}|+ n_k-j\}\cup \{l+1\})$ for $k>0$ 
and $t$ is the smallest $i$ such that $n_i=l+1$. 
We call $n(\Cal B)$ the {\it skeleton} of $\Cal B$. 

\proclaim{Lemma 1.3.4} Let $\Cal B$ be a string with skeleton $n(\Cal B)=
(n_1,\dots ,n_t).$  
Suppose that $1 \le j < t$,  $n_j \le i < n_{j+1}$, and 
$B_{i+1}\supseteq B_{i+2}, \dots , B_{n_{j+1}-1}$.  
Then $$n(T_i(\Cal B)) = (1,n_{j+1}-i,\dots ,n_t-i).$$
\endproclaim

\proclaim {Definition 1.3.5}
Let $B \subseteq \{1,\dots ,n\}$ and let $\Cal 
B=(B_1,\dots, B_l)$ be a string with skeleton
$n(\Cal B)=(n_1, \dots , n_t)$.
Define $d(B,\Cal B)$, $e(B,\Cal B)$ and $f(B, \Cal B)$, elements 
of $\{0,1,\dots ,n\}$, as follows:
\vskip 4 pt

$d(B,\Cal B) = e(B,\Cal B) =f(B, \Cal B)=0$ unless $B_1 = B 
\setminus \{b\}$ for some $b \in B$;
\vskip 4 pt

$d(B,\Cal B) = b$ if $B_1 = B \setminus \{b\}$ 
for some $b \in B;$
\vskip 4 pt

$e(B,\Cal B) = \max(B)$ if $d(B,\Cal B) \ne 0$;
\vskip 4 pt

$f(B,\Cal B) = e(B, \Cal B)$ if $d(B,\Cal B) \ne 0$ 
and either $t=2$ or else $t \ge 3$ and 
$|B_{n_2}| \ne |B| - n_2$;
\vskip 4 pt

$f(B,\Cal B) = \max(B \setminus (B \cap 
B_{n_2}))$ if $d(B,\Cal B) \ne 0, t \ge 3$ and 
$|B_{n_2}| = |B| - n_2.$
\endproclaim

\vskip 4 pt
\noindent {\bf Remark}. If $t=1$ (i.e., if $\Cal B =\emptyset$) then
$d(B, \Cal B)=0$ because there is no $B_1$.

\proclaim {Definition 1.3.6}
We say that $\Cal B$ is a {\it standard string} if
$B_{i+1}=B_i\setminus \{e(B_i, T_i(\Cal B))\}$ 
whenever $1\leq j < t$ and $n_j<i+1<n_{j+1}$.
\endproclaim

\proclaim {Definition 1.3.7}
We say that $\Cal B$ is a {\it reduced string} if 
$B_{i+1}=B_i\setminus \{f(B_i, T_i(\Cal B))\}$
whenever $1\leq j<t$ and $n_j<i+1< n_{j+1}$. 
\endproclaim

Let $Y=\{\Cal B \ |\ \Cal B \ {\text {is a standard string}}\}$
and
$Y^{\prime }=\{\Cal B \ |\ \Cal B \ {\text {is a reduced string}}\}$.

In the next section we will prove:

\proclaim {Theorem 1.3.8}
$r(Y)$ is a basis for $Q_n$.
\endproclaim

We will also show that
\proclaim {Theorem 1.3.9}
$r(Y^{\prime })$ is a basis for $Q_n$. 
\endproclaim

\noindent {\bf Remark}. Although the definition of $Y$ is more
transparent, it is easier
to work with $Y^{\prime}$. This is because the element $f(B, {\Cal B})$ 
depends on more detailed information about
the pair $(B, {\Cal B})$ than the element $e(B, {\Cal B}).$  
Thus $e(B,{\Cal B})$ depends only on $B$ and $B_1$, 
while $f(B,{\Cal B})$ depends on $B, B_1$ and $B_{n_2}.$  
(For example, $e(\{1,2,3\}, \{1,2\}, \{3\}) = e(\{1,2,3\},\{1,2\},\{1\}) = 3$ 
while $f(\{1,2,3\},\{1,2\},\{3\})$ $ = 2$ and $f(\{1,2,3\},\{1,2\},\{1\}) = 3$.)
This extra information carried by $f(B,{\Cal B})$ will be crucial for the 
induction arguments necessary to prove the linear independence of $Y^{\prime}.$

We first show that each of the sets
$r(Y)$ and $r(Y^{\prime })$ spans $Q_n$. 
We will observe that 
$r(Y)$ is linearly independent if and only if 
$r(Y^{\prime })$ is linearly
independent. We will then show that $r(Y^{\prime })$
is linearly independent and so
$r(Y)$ and $r(Y^{\prime })$ are bases for $Q_n$.
 
\smallskip

Let $A\subseteq \{1,\dots ,n\}$, $|A|=u$. Assume
$A=\{a_1,\dots , a_u\}$ where \break $a_1<a_2<\dots <a_u$.
Recall that,
for $1\leq j <u,$  $\Cal B(A:j)$ is the string \break
$(A, A\setminus \{a_u\}, \dots , A\setminus \{a_{u-j+2},
\dots , a_u\})$. Then $\Cal B(A:j)\in Y$ and any element of
$Y$ may be written as a juxtaposition of
$\Cal B(A_1:j_1), \dots , \Cal B(A_s:j_s)$ where
$A_1,\dots , A_s\subseteq \{1,\dots , n\}$,
$j_i\leq |A_i|$ for all $i$, and, for all $2\leq i \leq s$, we have either
$|A_i|\neq |A_{i-1}|-j_{i-1}$ or $A_i\nsubseteq A_{i-1}$.
Writing $r(A:j)=r(\Cal B(A:j))$ we see that Theorem 1.3.3 follows from Theorem 1.3.8.

\medskip

\subhead {1.4 Proof of the basis theorem}
\endsubhead

We begin by proving spanning results for $Y$ and $Y^{\prime}.$

\proclaim{Proposition 1.4.1} 
i) $r(Y\cap Q_{n,i})$
spans $Q_{n,i}$ for each $i\geq 0$.

ii) $r(Y^{\prime }\cap Q_{n,i})$
spans $Q_{n,i}$ for each $i\geq 0$.

\endproclaim

Our proof, as well as our subsequent arguments for linear independence, 
will depend on two partial orderings, $<$ and $<^{\prime}$ of the set of strings.

Let $\Cal B=(B_1,\dots ,B_l)$
and $n(\Cal B)=(n_1,\dots , n_t)$.
For each $1\leq q < k\leq l$ define $v_{k,q}(\Cal B)$,
usually written as $v_{k,q}$, to be 
$|B_k\setminus (B_q\cap B_k)|$.
For each 
$j$, $1\leq j <t$, and each $k$,
$n_j\leq k < n_{j+1}$, define a  $(k-n_j)$-tuple 
$v_k(\Cal B)$, usually written $v_k$, by
$v_k(\Cal B)=(v_{k,n_j}(\Cal B), v_{k,n_j+1}(\Cal B),
\dots , v_{k,k-1}(\Cal B))$.
Note that $v_{n_j}(\Cal B)$, $1\leq j \leq t-1$, is the empty sequence.
Set $v(\Cal B)=(v_1, v_2,\dots , v_l)$, the
juxtaposition of the sequences 
$v_1, v_2, \dots, v_l$.

Again, let $n(\Cal B)=(n_1,\dots , n_t)$. For every
$j$, $1\leq j <t$, and every $k$,
$n_j< k < n_{j+1}$, define $v_k^{\prime }(\Cal B)=v_k(\Cal B)$.
Let $v_1^{\prime }(\Cal B)=\emptyset $, the empty sequence, and,
for $1\leq j<t-1$, let $v_{n_{j+1}}^{\prime }(\Cal B)$,
usually written $v_{n_{j+1}}^{\prime }$, denote the 
$(n_{j+1}-n_j)$-tuple 
$(v_{n_{j+1},n_j}(\Cal B), v_{n_{j+1},n_j+1}(\Cal B),
\dots , v_{n_{j+1},n_{j+1}-1}(\Cal B))$.
Set $v^{\prime }(\Cal B)=(v_1^{\prime },v_2^{\prime },\dots , 
v_l^{\prime })$, the
juxtaposition of sequences $v_1^{\prime }, 
v_2^{\prime }, \dots, v_l^{\prime }$. 

We define two partial orderings, $<$ and $<^{\prime }$ on the set of all
strings
of a given length. This will be done in four steps.

\smallskip

{\it Step 1}. If $\Cal B=(B_1,\dots, B_l)$, $\Cal C=(C_1,\dots, C_l)$
and $|{\Cal B}|<|{\Cal C}|$, we say that $\Cal B < \Cal C$ 
and $\Cal B <^{\prime }\Cal
C$.

\smallskip

{\it Step 2}. Suppose $|\Cal B|=|\Cal C|$ and
$n(\Cal B)=(n_1, \dots , n_t)\neq n(\Cal C)=(m_1,\dots , m_s)$.  
Then, interchanging $\Cal B$ and $\Cal C$ 
if necessary, we may find some $j$, 
$1\leq j \leq \min \{s,t\}$, such that 
$n_1=m_1,\dots , n_{j-1}=m_{j-1}, n_j>m_j$ (i.e., $n({\Cal B}) > n({\Cal C})$ 
in the lexicographic order). 
In this case we say $\Cal B < \Cal C$ and $\Cal B <^{\prime }\Cal C$.

\smallskip

{\it Step 3}. Assume $|\Cal B|=|\Cal C|$ and $n(\Cal B)=n(\Cal C)$.
If $v(\Cal B) < v(\Cal C)$ in the lexicographic order, then
$\Cal B < \Cal C$. If $v^{\prime }(\Cal B) < v^{\prime }(\Cal C)$,
then $\Cal B <^{\prime }\Cal C$.
\smallskip

{\it Step 4}. Let $||B||=\sum _{i\in B}i$ and
$||\Cal B||=(||B_1||, ||B_2||, \dots ||B_l||)$.
Assume $|\Cal B|=|\Cal C|$ and  $n(\Cal B)=n(\Cal C)$.
If $v(\Cal B)=v(\Cal C)$, and $||\Cal B||<||\Cal C||$,
then $\Cal B < \Cal C$. If $v^{\prime }(\Cal B) =
v^{\prime }(\Cal C)$ and $||\Cal B||<||\Cal C||$,
then $\Cal B <^{\prime}\Cal C$.

\smallskip




Proposition 1.4.1 follows from the following lemma. 

\proclaim{Lemma 1.4.2} Let $\Cal B=(B_1,\dots , B_l)$.

(a) If $\Cal B \notin Y$, then $r(\Cal B)$ is a linear
combination of monomials $r(\Cal C)$ with $\Cal C < \Cal B$.

(b) If $\Cal B \notin Y^{\prime }$, then $r(\Cal B)$ is a linear
combination of monomials $r(\Cal C)$ with 
$\Cal C <^{\prime }\Cal B$.

\endproclaim

\noindent {\bf Remark}.  This lemma shows that the partial 
ordering $<$ is closely associated 
with the set $Y$ which is defined using the elements $e(B,{\Cal B})$, 
and that the partial ordering $<^{\prime}$ is closely associated 
with the set $Y^{\prime}$ which is defined using the elements $f(B,{\Cal B}).$
Thus (see the remark following the statement of Theorem 1.3.9) $<^{\prime}$ 
reflects more detailed structure of strings than $<$.  Because of this, 
the inductive arguments which prove independence are based on $<^{\prime}.$

{\it Proof of the lemma}. We will prove only part (a).
Part (b) can be proved similarly.
First, suppose $v(\Cal B)\neq
(0, \dots , 0)$. We will show that $r(\Cal B)$ is a
linear combination of monomials $r(\Cal C)$ 
with $\Cal C < \Cal B$.

Since $v(\Cal B)\neq (0, \dots , 0)$ we may find some $j$,
$1\leq j < t$, and some $k$,
$n_j < k < n_{j+1}$, such  that 
$v_k\neq (0,\dots , 0)$ and  
$v_m = (0,\dots , 0)$ or $\emptyset $
whenever $m < k$. Since $v_k\neq (0, \dots , 0)$
we may find some $s$, $n_j\leq s < k$, such that
$v_{k,s}\neq 0$ and $v_{k,q}=0$ whenever
$n_j\leq q <s$.

Now if $n_j \leq k_1 < k_2 < k$, we have $v_{k_2}=
(0, \dots , 0)$, and so $v_{k_2, k_1}=0$.  Thus
$B_{k_1}\supseteq B_{k_2}$, and so one has
$B_{n_j}\supseteq B_{n_j+1}\supseteq \dots \supseteq 
B_{k-1}$. By the definition of $n_{j+1}$, it follows that
$B_{n_j+i+1}=B_{n_j+i} \setminus \{b_i\}$, where
$b_i\in B_{n_j+i}$ for $i=0, \dots , k-n_j-2$.

Note that, for $n_j\leq q \leq k-2$, 
$v_{k, q+1}-v_{k, q}=0$ if $b_{q-n_j}\notin B_k$ and
$v_{k, q+1}-v_{k, q}=1$ if $b_{q-n_j}\in B_k$.
Thus $0\leq v_{k, q+1} - v_{k, q}\leq 1$.

Now suppose 
$v_k=(0, \dots , 0, v_{k,s},\dots , v_{k, k-1})=
(0, \dots , 0, 1,2,\dots , k-s)$, i.e., 
$v_{k, s-1+i}=i$ for $1\leq i \leq k-s$.
Then we have $b_{s-1-n_j}, \dots , b_{k-2-n_j}\in B_k$.
Since $|B_k|=|B_{k-1}|-1 < |B_{k-2}|$,
we can find $c\in B_{k-2}$, $c\notin B_k$. Set
$B_{k-1}^{\prime}=B_{k-2}\setminus \{c\}$.

Let $\Cal B^{\prime}=(B_1, \dots , B_{k-2},
B_{k-1}^{\prime}, B_k, \dots B_l)$. Then since
$|B_{k-1}^{\prime }|=|B_{k-1}|$, we have
$|\Cal B^{\prime}|=|\Cal B|$. Also since
$B_{k-1}^{\prime }\subseteq B_{n_j}$, we have
 $n(\Cal B^{\prime})=n(\Cal B)$. Clearly $v_m(\Cal B) = v_m(\Cal B')$ whenever $m < k-1.$
Since $B_{k-1}^{\prime }\subseteq B_{k-2} \subseteq 
\dots \subseteq B_{n_j}$, $v_{k-1}(\Cal B^{\prime })
=v_{k-1}(\Cal B_k)$ is either $(0, \dots , 0)$ or
$\emptyset $. 

We complete this case by noting that 
$v_{k,h}(\Cal B^{\prime })=v_{k,h}(\Cal B)$ if
$h < k-1$ and that $v_{k,k-1}(\Cal B^{\prime })=
v_{k, k-1}(\Cal B)-1$. Thus $v(\Cal B^{\prime })
< v(\Cal B)$ and so $\Cal B^{\prime } < \Cal B$. 
Since the quadratic relations (1.2) show that
$r(\Cal B^{\prime})$ equals $r(\Cal B)$ modulo
terms of lower degree, the lemma is proved in this case.

Next suppose there is some $p$, $1\leq p <k-s$, so that
$v_{k, s+i-1} \ne i$ for $1\leq i\leq p$ and
$v_{k, s+p}=v_{k, s+p-1}$. Hence 
$b_{s-1-n_j}, \dots, b_{s+p-2-n_j}\in B_k$, 
$b_{s+p-1-n_j}\notin B_k$. Set 
$B_{s+p-1}^{\prime }=B_{s+p-2}\setminus \{b_{s+p-1-n_j}\}$
and $(B_1, \dots , B_{s+p-2},
B_{s+p-1}^{\prime}, B_{s+p}, \dots B_l) \break = \Cal B^{\prime}$.

As before $|B_{s+p-1}^{\prime}|=|B_{s+p-1}|$,
so $|\Cal B^{\prime }|=|\Cal B|$
and $B_{s+p-2}\supseteq B_{s+p-1}^{\prime }
\supseteq B_{s+p}$ so $n(\Cal B^{\prime})=n(\Cal B)$.
Clearly
$v_u(\Cal B^{\prime})=v_u(\Cal B)$ whenever
$u<s+p-1$. Furthermore, since
$B_{s+p-2}\supseteq B_{s+p-1}^{\prime }\supseteq B_{s+p}$,
 we have $v_u(\Cal B^{\prime})=v_u(\Cal B)=(0, \dots , 0)$
whenever $s+p-1\leq u \leq k-1$.

To complete this case, note that $v_{k,h}(\Cal B^{\prime})=v_{k,h}(\Cal B)$ if
$h<s+p-1$ and that $v_{k,s+p-1}(\Cal B^{\prime})=
v_{k, s+p-1}(\Cal B)-1$. Thus, $v(\Cal B^{\prime}) < v(\Cal B)$
and so $\Cal B^{\prime} <\Cal B$.  As before, the quadratic relations (1.2) show that
$r(\Cal B^{\prime})$ equals $r(\Cal B)$ modulo
terms of lower degree, and so the lemma is proved in this case.

Suppose finally that $v(\Cal B)=(0, \dots , 0)$. This means
that $B_1\supseteq B_2\supseteq \dots \supseteq B_{n_2-1}$,
$B_{n_2}\supseteq B_{n_2+1}\supseteq \dots \supseteq B_{n_3-1},
\dots $.

Then $B_{i+1}=B_i\setminus \{d(B_i, T_i(\Cal B))\}$
whenever $1\leq j <t$ and $n_j<i+1<n_{j+1}$.
Since $\Cal B\notin Y$ we have a pair $(j,i)$,
$1\leq j <t$, $n_j<i+1<n_{j+1}$, such that
$d(B_i, T_i(\Cal B)) < e(B_i, T_i(\Cal B))$.
Assume that $j$ is minimal and that, for this fixed value of $j$,
$i$ is maximal. Set $B_{i+1}^{\prime }=B_i\setminus
\{e(B_i, T_i(\Cal B))\}$ and
$\Cal B^{\prime }=(B_1, \dots , B_i,
B_{i+1}^{\prime }, B_{i+2}, \dots , B_l)$. Since
$|B_{i+1}^{\prime}|=|B_{i+1}|$, 
we have $|\Cal B^{\prime}|=|\Cal B|$. Since
$B_{i+1}^{\prime}\subseteq B_{i}$ we have
$n(\Cal B^{\prime})=n(\Cal B)$.

By the maximality of $i$, either $i+1=n_{j+1}$ or 
$B_{i+2}=B_{i+1}\setminus \{e(B_{i+1}, T_{i+1}(\Cal B))\}$.
In the latter case, since $e(B_i, T_i(\Cal B))
\neq d(B_i, T_i(\Cal B))$, we have
$e(B_{i+1}, T_{i+1}(\Cal B))=e(B_i, T_i(\Cal B))$ and so
$B_{i+2}\subseteq B_{i+1}^{\prime }$. Thus, in either case,
$v(\Cal B^{\prime})=v(\Cal B)=(0, \dots , 0)$. 

Finally $||B_{i+1}^{\prime }||=||B_{i+1}|| +
d(B_i, T_i(\Cal B)) - e(B_i, T_i(\Cal B)) <
||B_{i+1}||$ so
$||\Cal B^{\prime }||<||\Cal B||$ and
$\Cal B^{\prime }<\Cal B$. As before, the application
of the quadratic relations (1.2) completes the proof of this case and of the lemma.

The following result is obvious.

\proclaim{Lemma 1.4.3} Let $n=(n_1,\dots , n_t)$ where
$1=n_1<n_2< \dots <n_t=l+1$ and let 
$B_{n_1},\dots , B_{n_{t-1}}\subseteq \{1,\dots , n\}$
be such that for each $i$, $1\leq i <t-1$, either
$B_{n_{i+1}}|\nsubseteq B_{n_i}$ or
$|B_{n_{i+1}}|\neq |B_{n_i}| -n_i+n_{i+1}$. Then there is
a unique $\Cal C=(C_1, \dots , C_l)\in Y$ such that
$n(\Cal C)=n$ and $C_{n_i}=B_{n_i}$ for
$1\leq i\le t$ and a unique 
$\Cal C^{\prime }=(C_1^{\prime }, \dots , C_l^{\prime })\in Y^{\prime }$
such that $n(\Cal C^{\prime })=n$ and $C_{n_i}^{\prime }=B_{n_i}$
for $1\leq i <t$. Furthermore, every $\Cal C\in Y$ and
every $\Cal C^{\prime } \in Y^{\prime }$ occurs in this way.
\endproclaim

\proclaim{Corollary 1.4.4} For all $i\geq 0,$
$|r(Y)\cap Q_{n,i}|=|r(Y^{\prime })\cap Q_{n,i}|$.
\endproclaim

{\it Proof}. Each of these is in one-to-one
correspondence with the same set of sequences of integers and
subsets of $\{1,\dots , n\}$.

\proclaim{Corollary 1.4.5} $r(Y)$ is linearly independent if and
only if $r(Y^{\prime })$ is linearly independent.
\endproclaim

We will presently prove Theorem 1.3.9.
Corollary 1.4.5 shows that this theorem is equivalent to
Theorem 1.3.8.

We will need the following technical lemma.

\proclaim{Lemma 1.4.6} Assume $B \subseteq \{1,\dots,n\}, {\Cal B} = (B_1,....,B_l) \in Y'$ 
and $d(B, {\Cal B}) \ne f(B,{\Cal B}).$  
Set ${\Cal E} = (B \setminus \{f(B, {\Cal B})\}, B_2,....,B_l).$  
Then $ {\Cal E} <' {\Cal B}.$
\endproclaim

{\it Proof}.
As usual, write $n({\Cal B}) = (n_1,\dots,n_t).$  We claim that $n({\Cal E}) \ge n({\Cal B})$ 
and that if $n({\Cal E}) = n({\Cal B})$ then $v'({\Cal E}) \le v'({\Cal B}).$  
Since $T_1({\Cal E}) = T_1({\Cal B})$, this is immediate if $n_2 = 2.$  Thus 
we may assume $n_2  > 2.$  It is then enough to show that 
$B \setminus \{f(B,{\Cal B})\} \supseteq B_2.$  
Since $B_1 = B \setminus \{d(B,{\Cal B})\}$ 
and $d(B, {\Cal B}) \ne f(B,{\Cal B})$, we have $f(B,{\Cal B}) \in B_1.$ 
The maximality of $f(B,{\Cal B})$ then shows that 
$f(B,{\Cal B}) = f(B_1,T_1({\Cal B})),$ proving the claim.

Now suppose $n({\Cal E}) = n({\Cal B})$ and $v'({\Cal E}) = v'({\Cal B}).$  
Note that if $t \ge 3$ and $|B_{n_2}| = |B| - n_2$, we have $d(B,{\Cal B}) \notin B \cap B_{n_2}$.  
Then, in any case, the maximality of $f(B,{\Cal B})$ implies $f(B,{\Cal B}) > d(B,{\Cal B})$ 
and so $||{\Cal E}|| < ||{\Cal B}||$, proving the lemma.

We now begin the proof of Theorem 1.3.8. This is similar
to the standard proof of the Poincar\'e-Birkhoff-Witt Theorem
(see e.g., \cite {J}).

Let $F$ be the free algebra on generators $s_B$, 
$ B\subseteq \{1, \dots , n\}$, and let $V^{\prime}$ be a vector space
with basis $Y^{\prime}$. We will inductively
define an action of $F$ on $V^{\prime }$ by using the
following lemma.
\proclaim {Lemma 1.4.7} There is a unique 
action of $F$ on $V^{\prime }$ 
$$(s,v)\mapsto s*v$$
for $s\in F$, $v\in V^{\prime }$, 
such that $s_{\emptyset} = 0$ and for $\emptyset \ne B \subseteq \{1,\dots,n\}$ and ${\Cal B} = (B_1,\dots,B_l) \in Y'$ we have:

i) if
$d(B, \Cal B) = f(B, \Cal B)$, then $s_B*(\Cal B)=(B, B_1, \dots ,B_l)$;

\noindent and

ii) if $d(B, \Cal B)\neq f(B, \Cal B)$, then $$ s_B*\Cal B=
s_B*s_{B\setminus \{f(B, \Cal B)\}}*T_1({\Cal B})$$
$$- (s_{B\setminus \{d(B, \Cal B)\}}-s_{B\setminus \{f(B, \Cal B)\}})
*s_{B\setminus \{d(B, \Cal B), f(B, \Cal B)\}}*T_1({\Cal B})$$
$$+ s_{B\setminus \{d(B, \Cal B)\}}*s_{B\setminus \{d(B, \Cal B)\}}
*T_1({\Cal B})$$
$$- s_{B\setminus \{f(B, \Cal B)\}}*s_{B\setminus \{f(B, \Cal B)\}})
*T_1({\Cal B})$$.

\noindent In fact,

iii) if $d(B,{\Cal B}) \ne f(B,{\Cal B})$, then $s_B*{\Cal B}$ is a linear combination of strings ${\Cal D} \in Y'$ 
such that $l({\Cal D}) = 1 + l({\Cal B})$, $|{\Cal D}| \le |B| + |{\Cal B}|$, and ${\Cal D} < (B,B_1,\dots,B_l).$ 
\endproclaim

{\it Proof}. 
Partially order the set of pairs $(B, \Cal B)$, 
$\emptyset \neq B\subseteq \{1, \dots , n\}$,
$\Cal B \in Y^{\prime }$, by $(B, \Cal B)<(C, \Cal C)$ if
$|B| + |{\Cal B}| < |C| + |{\Cal C}|$ or if 
$|B| + |{\Cal B}| = |C| + |{\Cal C}|$ and ${\Cal B} <' {\Cal C}$.
We will define $s_B*\Cal B$ and prove that (i)-(iii)
hold inductively.

Define $s_B*{\emptyset} = (B)$ and $s_{\emptyset}*{\Cal B} = 0$ 
for all $\emptyset \ne B \subseteq \{1,\dots,n\}$ and all ${\Cal B} \in Y'$.
Now assume that for some pair $(C, {\Cal C})$, $s_B*{\Cal B}$ has been defined
for all pairs $(B, {\Cal B}) < (C,{\Cal C})$ and that (i)-(iii) hold for all such pairs.  
We will define $s_C*{\Cal C}$ and show that (i) - (iii) are satisfied for $s_C*{\Cal C}$. 

If $d(C,{\Cal C}) = f(C,{\Cal C})$ we define $s_C*({\Cal C}) = (C,C_1,\dots,C_l).$  
If $d(C,{\Cal C}) \ne f(C,{\Cal C})$, then, by 
Lemma 1.4.6, $(C \setminus \{f(C,{\Cal C})\},C_2,\dots,C_l) <' {\Cal C}.$
Furthermore, \break
$(C \setminus \{f(C,{\Cal C})\},T_1({\Cal C})) < (C,{\Cal C})$ 
and so, by induction (using (iii)), we have that 
$s_{C \setminus \{f(C,{\Cal C})\}}*T_1({\Cal C})$ is a linear 
combination of strings ${\Cal D} \in Y'$ such that
$l({\Cal D}) = 1 + l({\Cal B})$ and 
${\Cal D} <' (C \setminus \{f(C,{\Cal C})\},C_2,\dots,C_l).$  
Thus  $s_{C \setminus \{f(C,{\Cal C})\}}*T_1({\Cal C})$ is a 
linear combination of strings ${\Cal D} <' {\Cal C}$ and so
all the terms occuring in the expression (ii) for $s_C*{\Cal C}$ 
are defined by induction.  We use this expression to define  
$s_C*{\Cal C}$.

It remains to show that $(C,{\Cal C})$ satisfies (iii).  
Since $d(C,{\Cal C}) \ne f(C,{\Cal C})$ then, using (ii) and 
the induction assumption, we see that it is sufficient to show that 
$s_C*s_{C \setminus \{f(C,{\Cal C})\}}*T_1({\Cal C})$ is a 
linear combination of strings $\Cal D$ satisfying 
the conditions of (iii).  Since, as noted above, 
$s_{C \setminus \{f(C,{\Cal C})\}}*T_1({\Cal C})$ is a linear combination of
strings $<^{\prime} {\Cal C}$, the result follows by induction and
the lemma is proved.

\medskip

The map $s_B \mapsto r(B)$ extends to an epimorphism of $F$
onto $Q_n$.  We will show (Lemma 1.4.9) that the quadratic relations (1.2) are in 
the kernel of this epimorphism. Hence the action of $F$ on $V^{\prime}$ induces 
an action of $Q_n$ on $V^{\prime }$. If ${\Cal B} \in Y'$, then, by Lemma 1.4.7(i), we
have ${\Cal B}*{\emptyset} = {\Cal B}$.  Since $Y'$ is, by definition, 
a linearly independent subset of $V'$, $r(Y')$ is linearly independent and 
Theorem 1.3.9 follows.
 
For $B \subseteq \{1,\dots,n\}$ 
and $i,j \in B$ define $\phi(B,i,j) \in F$ by 
$$\phi(B,i,j) = - (s_{B\setminus i}-s_{B\setminus j})s_{B\setminus \{i,j\}} 
+ s_{B\setminus i}^2 - s_{B\setminus j}^2.$$

\proclaim {Lemma 1.4.8} 

(a) $\phi (B,i,j)=-\phi (B,j,i)$;

(b) $\phi (B\setminus i,j,k) + \phi (B\setminus j,k,i) +
\phi (B\setminus k,i,j)=0$;

(c) $\phi (B,i,j) + \phi (B,j,k) + \phi (B,k,i)$ \hfill \break $  =
-(s_{B\setminus i} - s_{\setminus j})s_{B\setminus \{i,j\}}
- (s_{B\setminus j} - s_{\setminus k})s_{B\setminus \{j,k\}}
- (s_{B\setminus k} - s_{\setminus i})s_{B\setminus \{i,k\}}. $
\endproclaim

{\it Proof}. These are all immediate from the definition.

\proclaim {Lemma 1.4.9}
$$
(s_C(s_{C\setminus i}-s_{C\setminus j})
-\phi(C,i,j))*\Cal C =0
$$
for all strings $\Cal C \in Y'$ and for all $i,j \in C \subseteq \{1,\dots,n\}.$
\endproclaim

Assume ${\Cal C} = (C_1,\dots,C_l).$  We prove the lemma by induction 
on the greater (with respect to $<^{\prime}$) of the two strings 
$(C, C\setminus i, C_1,...,C_l)$ and 
$(C, C\setminus j, C_1,...,C_l)$. 
Note that the lemma holds whenever $|C| + |\Cal C| \le 1$.  Assume that for 
some $u,v \in B \subseteq \{1,\dots,n\}$ and some ${\Cal B} = (B_1,\dots,B_s) \in Y'$, 
the lemma holds whenever both $(C, C\setminus i, C_1,...,C_l)$ and 
$(C, C\setminus j, C_1,...,C_l)$ are less than the greater of the 
two strings $(B, B\setminus u, B_1,...,B_s)$ and $(B, B\setminus v, B_1,...,B_s)$. We will show that 
$(s_B(s_{B\setminus u}-s_{B\setminus v})*{\Cal B} 
= \phi(B,u,v))*\Cal B.$

We will need two preliminary steps.

{\it Step 1}. Let $u,v\in B\subseteq \{1,\dots , n\}$ and $\Cal B$
be as above. Let $w\in B$. Then 
$$
(\phi(B,u,v) + \phi(B,v,w) + \phi(B,w,u))*\Cal B =0.
$$

{\it Proof of Step 1}. Using Lemma 1.4.8(c) we see that
$(\phi(B,u,v) + \phi(B,v,w) + \phi(B,w,u))*\Cal B =
-(s_{B\setminus u}(s_{B\setminus \{u,v\}} - s_{B\setminus \{u,w\}}) +
s_{B\setminus v}(s_{B\setminus \{v,w\}} - s_{B\setminus \{u,v\}}) +
s_{B\setminus w}(s_{B\setminus \{u,w\}} - s_{B\setminus \{v,w\}}))*\Cal B$.

By induction, this is equal to
$-(\phi(B\setminus u,v,w) + \phi(B\setminus v,w,u) + \phi(B\setminus w,u,v))*\Cal B$.
This is equal to $0$ by Lemma 1.4.8(b), and so Step 1 is complete.

{\it Step 2}.  Let $u,v\in B\subseteq \{1,\dots , n\}$ and $\Cal B$
be as above. Assume $|B_1| = |B| - 2.$ Let $w\in B$. Then 
$$
(\phi(B,u,v)s_{B\setminus \{u,v\}} + \phi(B,v,w)s_{B\setminus \{v,w\}}
 + \phi(B,w,u)s_{B\setminus \{u,w\}})*(B_2,\dots , B_l) =0.
$$

{\it Proof of Step 2}. Using the definition of $\phi $ we may rewrite the
left-hand side as
$$
(-(s_{B\setminus u}-s_{B\setminus v})s_{B\setminus \{u,v\}}^2
- (s_{B\setminus v}-s_{B\setminus w})s_{B\setminus \{v,w\}}^2
- (s_{B\setminus w}-s_{B\setminus u})s_{B\setminus \{u,w\}}^2
$$ $$ -s_{B\setminus u}^2s_{B\setminus \{u,v\}}
+ s_{B\setminus v}^2s_{B\setminus \{u,v\}}
- s_{B\setminus v}^2s_{B\setminus \{v,w\}}
$$ $$ +s_{B\setminus w}^2s_{B\setminus \{v,w\}}
- s_{B\setminus w}^2s_{B\setminus \{u,w\}}
+ s_{B\setminus u}^2s_{B\setminus \{u,w\}})*T_1({\Cal B}).
$$

This is equal to 
$$
(s_{B\setminus u}(s_{B\setminus u}(s_{B\setminus \{u,v\}}
-s_{B\setminus \{u,w\}}) - s_{B\setminus \{u,v\}}^2 +
s_{B\setminus \{u,w\}}^2) $$ $$+
s_{B\setminus v}(s_{B\setminus v}(s_{B\setminus \{v,w\}}
-s_{B\setminus \{u,v\}}) - s_{B\setminus \{v,w\}}^2 +
s_{B\setminus \{u,v\}}^2) $$ $$+
s_{B\setminus w}(s_{B\setminus w}(s_{B\setminus \{u,w\}}
-s_{B\setminus \{v,w\}}) - s_{B\setminus \{u,w\}}^2 +
s_{B\setminus \{v,w\}}^2))*T_1({\Cal B}).
$$

By the induction assumption this is equal to
$$
(s_{B\setminus u}(\phi (B\setminus u,v,w)-s_{B\setminus \{u,v\}}^2
+s_{B\setminus \{u,w\}}^2) $$ $$+
s_{B\setminus v}(\phi (B\setminus v,w,u)-s_{B\setminus \{v,w\}}^2
+s_{B\setminus \{u,v\}}^2)  $$ $$+
s_{B\setminus w}(\phi (B\setminus w,u,v)-s_{B\setminus \{u,w\}}^2
+s_{B\setminus \{v,w\}}^2))*T_1({\Cal B}).
$$

By the definition of $\phi $ this is equal to
$$
-(s_{B\setminus u}(s_{B\setminus \{u,v\}}- s_{B\setminus \{u,w\}})
s_{B\setminus \{u,v,w\}} $$ $$+
s_{B\setminus v}(s_{B\setminus \{v,w\}}- s_{B\setminus \{u,v\}})
s_{B\setminus \{u,v,w\}} $$ $$+
s_{B\setminus w}(s_{B\setminus \{u,w\}}- s_{B\setminus \{v,w\}})
s_{B\setminus \{u,v,w\}})*T_1({\Cal B}).
$$

By the induction assumption this is equal to
$$
- ((\phi (B\setminus u,v,w) + \phi (B\setminus v,w,u)
+ \phi (B\setminus w,u,v))s_{B\setminus \{u,v,w\}}*
T_1({\Cal B}).
$$

By Lemma 1.4.8(b) this is equal to $0$, and so Step 2 is complete.

We are now ready to prove the lemma. There are three
cases depending on the relations among $B$, $B_1$, $u$
and $v$.

{\it Case 1}. Either ${\Cal B} = \emptyset$, or $B_1\neq B\setminus \{a,b\}$ for any
$a,b\in B$, $a\neq b$, or $B_1= B\setminus \{a,b\}$
with $\{a,b\}\cap \{u,v\}=\emptyset $. Then, by the definition of *,
$$
s_B(s_{B\setminus u}-s_{B\setminus v})*{\Cal B}
$$
$$= s_B*(B \setminus u, B_1,\dots,B_s) - s_B*(B \setminus v, B_1,\dots,B_s)  $$
$$= s_B*s_{B \setminus \{f(B, (B \setminus u, B_1,\dots,B_s))\}}*{\Cal B} 
+ \phi(B,u,f(B,(B \setminus u, B_1,\dots,B_s))* {\Cal B}  $$ 
$$- s_B*s_{B \setminus \{f(B, (B \setminus v, B_1,\dots,B_s))\}}*{\Cal B} 
- \phi(B,v,f(B,(B \setminus v, B_1,\dots,B_s))* {\Cal B}.$$ 
However, $$f(B,(B \setminus u, B_1,\dots,B_s)) = 
f(B, (B \setminus v,B_1,\dots,B_s))$$ since, in either case, this is the 
largest element of $B$ if ${\Cal B} = \emptyset$ or if $|B_1| \ne |B| - 2$, 
while if $|B_1| = |B| - 2$ it is the largest element 
of $(B \setminus (B \cap B_1))$ if $|B_1| = |B| - 2$.  Thus, setting 
$w = f(B,(B \setminus u, B_1,\dots,B_s)),$ our expression becomes 
$$(\phi(B,u,w) - \phi(B,v,w))*{\Cal B}.$$  By Step 1 this is equal 
to $\phi(B,u,v)*{\Cal B},$ as required.
\smallskip
{\it Case 2}.  $B_1 = B \setminus \{a,u \}$ where $a \in B, a \ne u, a \ne v.$  
Then, by the definition of $*$, 
$$s_B*(s_{B \setminus u} - s_{B \setminus v})*{\Cal B}  $$
$$= s_Bs_{B \setminus u} * {\Cal B} - s_B*s_{B \setminus 
f(B,(B \setminus v, B_1,\dots,B_s))}*{\Cal B} $$ $$
- \phi(B,v,f(B,(B \setminus v, B_1,\dots,B_s)))*{\Cal B}.$$
Note that $f(B,(B \setminus v,B_1,\dots,B_s))$ is the largest 
element of $B \setminus (B \cap B_1).$  Thus $f(B, (B \setminus v.B_1,\dots,B_s)) \ne v$
and $f(B, (B \setminus v,B_1,\dots,B_s)) \ge u.$  
If \break $f(B,(B\setminus v,B_1,\dots,B_s)) = u$, then our expression is 
obviously equal to \break 
$-\phi(B,v,u)*{\Cal B} = \phi(B,u,v)*{\Cal B}$, 
as required.  Hence we may assume that 
$f(B, (B \setminus v, B_1,\dots,B_s)) = a > u.$ Now 
$$(B, B \setminus f(B, (B \setminus v, B_1,\dots,B_s),B_1,\dots,B_s)) =
(B, B\setminus v, B_1,\dots,B_s).$$  Furthermore, 
$(B, B\setminus u, B_1,\dots,B_s)<^{\prime} (B, B\setminus v, B_1,\dots,B_s)$ and  \break
$(B, B\setminus a, B_1,\dots,B_s)<^{\prime} (B, B\setminus v, B_1,\dots,B_s).$
Hence the induction assumption shows that
$$s_B * s_{B \setminus u} * {\Cal B} - s_B * 
s_{B \setminus\{f(B,(B \setminus v, B_1,\dots,B_s))\}}*{\Cal B} 
= \phi(B,u,a)*{\Cal B}$$ and 
so our expression is equal to $$\phi(B,u,a) 
- \phi(B, v, a) =  \phi(B,u,v),$$ as required.
\smallskip

{\it Case 3}.  $B_1 = B \setminus \{u,v \}$.  

If $t=2$ or if $t > 2$ and $|B_{n_2}| \ne |B_1| + 1 - n_2$, let $B_{n_2}'$ denote $\emptyset$.  Otherwise, let $B_{n_2}'$
denote $B_{n_2}.$  
Assume, without loss of generality, that $u > v$.
\smallskip
{\it Subcase 3a}.   Let $a$ denote the largest element of $B$ not 
contained in $B_{n_2}'$.  Assume $a \ne u,v.$  Then, by induction,
$$s_B s_{B \setminus u}*(B \setminus \{u,v\},B_2,\dots,B_l) -
s_B s_{B \setminus v}*(B \setminus \{u,v\},B_2,\dots,B_l) $$
$$= s_B s_{B \setminus u}s_{B \setminus \{ u,a \}} * T_1({ \Cal B }) + 
s_B \phi({B \setminus u},v,a) * T_1({ \Cal B })  $$
$$-s_B s_{B \setminus v}s_{B \setminus \{v,a\}}*T_1({\Cal B}) -
s_B \phi({B \setminus v},u,a) * T_1({\Cal B}).$$
By Lemma 1.4.7, this is equal to 
$$s_B s_{B \setminus a}s_{B \setminus \{ u,a \}} * T_1({ \Cal B }) + 
 \phi(B,u,a)s_{B \setminus \{u,a\}}*T_1({\Cal B}) +
s_B \phi({B \setminus u},v,a) *T_1( { \Cal B })  $$
$$- s_B s_{B \setminus a}s_{B \setminus \{v,a\}}*T_1({\Cal B}) -
\phi(B,v,a)s_{B \setminus \{v,a \}}*T_1({\Cal B})
- s_B \phi({B \setminus v},u,a) * T_1({\Cal B}).  $$
By induction and Lemma 1.4.8, this is equal to
$$s_B \phi(B \setminus a,u,v)* T_1({\Cal B}) + s_B \phi(B 
\setminus u,v,a)*T_1( {\Cal B}) + s_B \phi(B \setminus v,a,u)* T_1({\Cal B})
$$
$$+ \phi(B,u,a)s_{B\setminus\{u,a\}}*T_1({\Cal B}) 
+ \phi(B,a,v)s_{B\setminus\{v,a\}}*T_1({\Cal B}).$$
By Lemma 1.4.8  and Step 2 this is equal to 
$$-\phi(B,v,u)s_{B\setminus\{u,v\}}*T_1({\Cal B}) 
= \phi(B,u,v)*(B\setminus \{u,v \},B_2,\dots,B_l)= \phi(B,u,v),$$ as required.

\vskip 8 pt

{\it Subcase 3b}.   Assume $u$ is the largest element of $B$ 
not contained in $B_{n_2}'$.  Let $b$ denote the largest
element of $B \setminus u$ not contained in $B_{n_2}'.$
Then, by induction,
$$s_B s_{B \setminus u}*(B \setminus \{u,v\},B_2,\dots,B_l) -
s_B s_{B \setminus v}*(B \setminus \{u,v\},B_2,\dots,B_l) $$
$$= s_B s_{B \setminus u}s_{B \setminus \{ u,b \}} * T_1({ \Cal B }) + 
s_B \phi({B \setminus u},v,b) * T_1({ \Cal B })  $$
$$- s_B*(B \setminus v,B \setminus \{u,v\},B_2,\dots,B_l). $$
By Lemma 1.4.7, this is equal to 
$$(B, B \setminus u, B \setminus \{u,b \},B_2,\dots,B_l)   
+ s_B \phi({B \setminus u},v,b) * T_1({ \Cal B })  $$
$$- s_B s_{B \setminus u}*(B \setminus \{u,v\},B_2,\dots,B_l)
- \phi(B ,v,u)*(B \setminus \{u,v \},B_2,\dots,B_l \}.$$
By induction and Lemma 1.4.8, this is equal to
$$s_B \phi(B \setminus u,b,v)* T_1({\Cal B}) + s_B \phi(B 
\setminus u,v,b)* T_1({\Cal B}) + \phi(B,u,v)*{\Cal B} =
\phi(B,u,v)*{\Cal B},$$ as required.

\vskip 8 pt

{\it Subcase 3c}. Assume $v$ is the largest element of $B$ not 
contained in $B_{n_2}'$.  Note that this implies $u \in
B_{n_2}.$ Let $c$ denote the largest element of $B \setminus v$ 
not contained in $B_{n_2}'.$
Then, by Lemma 1.4.7 and induction,
$$s_B s_{B \setminus u}*(B \setminus \{u,v\},B_2,\dots,B_l) -
s_B s_{B \setminus v}*(B \setminus \{u,v\},B_2,\dots,B_l) $$
$$= s_B* (B \setminus u,B \setminus \{ u,v \},B_2,\dots,B_l)  - 
s_B s_{B \setminus v} s_{B \setminus  \{v,c \}}*T_1({\Cal B}) 
- s_B \phi(B \setminus v,u,c)*T_1({\Cal B})  $$
$$= s_B s_{B \setminus v}*{\Cal B} + \phi(B,u,v)*{\Cal B}  $$
$$- s_B s_{B \setminus v}*(B \setminus \{v,c\},B_2,\dots,B_l)
- s_B \phi(B \setminus v,u,c)*T_1({\Cal B}) = $$
$$s_B \phi(B \setminus v,u,c)* T_1({\Cal B}) + \phi(B,u,v)*{\Cal B} - s_B \phi(B \setminus v,u,c)* T_1({\Cal B})  =
\phi(B,u,v)*{\Cal B},$$ as required.

This completes Case 3 and so completes the proof of Lemma 1.4.9 and hence of Theorem 1.3.9.

\subhead 1.5. Symmetric functions and subalgebras
of quadratic algebra $Q_n$ 
\endsubhead

Recall that the symmetric group, $S_n$, acts
on $Q_n$ by 
$\sigma (z_{A,i})=z_{\sigma (A),\sigma (i)}$
for $\sigma\in S_n$,  $A\subset \{1,\dots,n\}$,
$i\in \{1,\dots,n\}$, $i\notin A$. 

We are going to construct a family of $S_n$-invariant
elements in $Q_n$. To do this we
define inductively a family of expressions $\Lambda _k(A)$
in $Q_n$, $A\subset \{1,\dots,n\}$, $k=0,1,\dots, $ such that

i) $\Lambda _0(A)=1$ for all $A$,

ii) $\Lambda _k(\emptyset )=0$ for all $k\geq 1$,

iii) $\Lambda _k(A\cup i)=
\Lambda _k(A)+(\Lambda _1(A\cup i)-\Lambda _1(A))
\Lambda _{k-1}(A)$ for
all
$k\geq 1, i\notin A$.

\smallskip
It is easy to see that $\Lambda _k(A)=0$ if
$|A|\leq k-1$.

\noindent {\bf Example}. If $A=\{i\}$, then $\Lambda _1(A)=z_i$.
If $A=\{i,j\}$, then, by relations (1.1), 
$\Lambda _1(A)=z_i+z_{i,j}=z_j+z_{j,i}$,
$\Lambda _2(A)=z_{i,j}z_i=z_{j,i}z_j$. Note also that
$\Lambda _1(A)=r(A)$ in the notation of Section 1.2.

Set  $y_i=z_{12\dots i-1,i}$. From the relations (1.1) we have: 
\proclaim {Proposition 1.5.1} (a) The elements $\Lambda _k(A)$
do not depend on the ordering of $A$. In other words, 
if $\sigma (A)=A$, $\sigma \in S_n$,
then $\sigma (\Lambda _k(A))=
\Lambda _k(A)$ for all $k$.

(b) $$\Lambda _k(A)=
\sum _{i_1, \dots , i_k\in A, i_1>i_2>\dots >i_k}\
y_{i_1}y_{i_2}\dots y_{i_k}.
$$        
\endproclaim

In particular, the expressions 
$\Lambda _k(\{1,\dots , n\})$ are invariant under the
action of the symmetric group $S_n$. They can be written as
$$
\Lambda _k(\{1,\dots, n\})=
\sum _{n\geq i_1>i_2>\dots >i_k\geq 1}
y_{i_1}y_{i_2}\dots y_{i_k}.
$$

We will often write $\Lambda _k$ instead of
$\Lambda _k(\{1,\dots , n\})$. The commutative
analogues of the $\Lambda _k$'s are the 
well-known elementary symmetric functions in $z_1,\dots, z_n$.

\medskip

There is a bijective correspondence
between orderings $I=(i_1,\dots , i_n)$ of $\{1,\dots , n\}$
and $\sigma \in S_n$ such that
$\sigma (1)=i_1$,\dots, $\sigma (n)=i_n$.
The ordering $I$ defines a subalgebra $Q_{n, \sigma }$ in
$Q_n$ generated by the $r(\{i_1, \dots , i_k\})$, $k=1,\dots ,n$. 
It is evident that
$Q_{n, \sigma }$ is in fact a differential subalgebra of $Q_n$.
Also, $\tau (Q_{n, \sigma })=Q_{n, \tau \sigma }$
for all $\tau \in S_n$.

It was proved in \cite {GR3, GR5} that the polynomial
$P(t)\in Q_n[t]$ defined by the formula
$$
P(t)=(t-y_{n, \sigma })\dots (t-y_{1, \sigma })=
t^n - \Lambda _1t^{n-1} +\dots +(-1)^n\Lambda _n
$$ 
does not depend on the
permutation $\sigma $ and that the coefficients of $P(t)$ are 
$S_n$-invariant.

\proclaim{Theorem 1.5.2} The following subalgebras of the 
algebra $Q_n$ coincide:

(i) The intersection of all subalgebras $Q_{n, \tau }$;

(ii) The subalgebra of all $S_n$-invariant elements of $Q_{n, \sigma }$
for any permutation $\sigma $;

(iii) The subalgebra generated by the coefficients of
$P(t)$. 
\endproclaim

\noindent {\bf Example}. Let $n=2$. Denote by $Q_2'$ the subalgebra
of $Q_2$ generated by $z_1, z_{1,2}$, and by $Q_2''$ the
subalgebra generated by $z_2, z_{2,1}$. The theorem says
that the intersection $Q_2'\cap Q_2''$ is generated by the
$S_2$-invariant elements $z_1+z_{1,2}$ and 
$z_{1,2}\cdot z_1$.

It is interesting to note that Theorem 1.5.2 
reflects the fact that $Q_n$ is noncommutative.
Thus, let $\phi :Q_n\to L$ be any homomorphism of
$Q_n$ into a commutative integral domain $L$ which is
a $k$-algebra and assume that the elements 
$\phi (z_i)=t_i$ are distinct in $L$. Then, by 
Proposition 1.1.2, $\phi (z_{A,i})=t_i$ for all
$i\notin A$. Hence $\phi (Q_{n, \tau })=\phi (Q_n)$ is
the subalgebra generated by the $t_i$'s and so is
$\bigcap _{\tau }\phi (Q_{n, \tau })$.

On the other hand, $\phi (\bigcap _{\tau }Q_{n, \tau })$ is the image 
under $\phi $ of the subalgebra generated by all
$\Lambda _k$'s. But 
$\phi (\Lambda _k)=e_k(t_1,\dots , t_n)$, the elementary
symmetric function in commuting variables. Thus
$\phi (\bigcap _{\tau }Q_{n, \tau })$ is generated by
$e_k(t_1,\dots , t_n)$, $k=1,\dots , n$.  

Theorem 1.5.2 follows from the more general Theorem 1.5.5.

\proclaim{Conjecture 1.5.3} Let $w\in Q_n$ and $\sigma $
be a permutation of $\{1,\dots , n\}$. Assume that 
$w$ is a rational expression in
elements of $Q_{n, \tau }$. Then $w\in Q_{n, \tau } \in S_n$.
\endproclaim

The conjecture implies that $w\in Q_n$ is $S_n$-invariant if
and only if $w$ can be rationally expressed via elements of
$Q_{n, \tau }$ for all permutations $\tau \in S_n$.

Let $S$ be a subset of the symmetric group $S_n$. Then we may construct as follows the 
unique finest partition of
$\{1,\dots , n\}$ into $S$-invariant intervals.
Let $1=i_1 < \dots < i_{t(S)+1}=n+1 $
be such that for each $j, 1 \le j \le t(S), 
I_j(S) = \{i_j,...,i_{j+1}-1\}$ is a minimal, $S$-invariant interval.  

Recall that $Q_n$ has an increasing filtration $Q_{n,0}\subseteq Q_{n, 1}\subseteq
\dots $.  The associated graded algebra is denoted 
$\text {gr}\ Q_n=\sum_{s \ge 0} (\text {gr}\ Q_n)_{[s]} $ where
$(\text {gr}\ Q_n)_{[s]}=Q_{n,s}/Q_{n, s-1}$. For an
element $a\in Q_{n,s}\setminus Q_{n, s-1}$ denote by $\bar a$ the
corresponding element in  $(\text {gr}\ Q_n)_{[s]}$.

Let $P$ denote the set of all sequences $((i_1, j_1),\dots ,
(i_s, j_s))$ of pairs of integers with $n\geq i_q\geq j_q\geq 1$
for $1\leq q\leq s$ and let $P^{\prime }$ denote the subset
of $P$ consisting of all such sequences with $i_{q+1}\neq
i_q-j_q$ for $1\leq q < s$. In the notation of Section 1.3 let
$
G((i_1,j_1),\dots , (i_s, j_s))$ denote the set of all products 
$\bar r(B_1:j_1) \dots \bar r(B_s, j_s)$ such that
$((i_1, j_1),\dots , (i_s, j_s))\in P$,
$|B_q|=i_q $ for $ 1\leq q\leq s$ and for each
$q$, $1\leq q< s$ either $i_{q+1}\neq i_q-j_q$
or $B_{q+1}\nsubseteq B_q$.

Clearly, we have
\proclaim{Lemma 1.5.4} (a) $\bar r(Y)=
\bigcup _{((i_1,j_1),\dots , (i_s, j_s))\in P}
G((i_1,j_1),\dots , (i_s, j_s))$.

(b) $G((i_1,j_1),\dots , (i_s, j_s))\cap
G((i_1',j_1'),\dots , (i_s', j_s'))=\emptyset $

unless $((i_1,j_1),\dots , (i_s, j_s))=
((i_1',j_1'),\dots , (i_s', j_s'))$.

(c) $S_n$ permutes the elements of 
$G((i_1,j_1),\dots , (i_s, j_s))$.
\endproclaim

\proclaim{Theorem 1.5.5} (a) Let $\tau \in S\subseteq S_n$.
Then $\tau (\Lambda _m (I_j(S)))=\Lambda _m(I_j(S))$
for any $j$, $1\leq j \leq t(S)$, $m\leq |I_j(S)|$.

(b) $\bigcap _{\tau \in S}Q_{n, \tau }$ is generated by
$\{\Lambda _m (I_j(S))\ |\ 1\leq j\leq t(S),
1\leq m\leq |I_j(S)|\}$.
\endproclaim

{\it Proof}. (a) Write $J_j(S)=I_1(S)\cup \dots \cup I_j(S)$.
By Proposition 1.5.1(a) $\tau (\Lambda _m(J_j(S)))=
\Lambda _m(J_j(S))$ for $1\leq m \leq |J_j(S)|$. In particular,
taking $j=1$, we have  $\tau (\Lambda _m(I_1(S)))=
\Lambda _m(I_1(S))$ for $1\leq m \leq |I_1(S)|$.

Clearly $$
\Lambda _m(J_j(S))=\Lambda _m(I_j(S)\cup J_{j-1}(S))=\sum _{k=0}^m
\Lambda _k(I_j(S))\Lambda _{m-k}(J_{j-1}(S))
,$$ 
and so
$$
\Lambda _m(I_j(S))=\Lambda _m(J_j(S)) - \sum _{k=0}^{m-1}
\Lambda _k(I_j(S))\Lambda _{m-k}(J_{j-1}(S)).
$$ 
The result now follows by double induction, first on $j$ and
then on $m$.

(b) This proof is based on Lemma 1.5.4. Note that
$Q_{n, \sigma }$ is generated by $r(\sigma \{1\})$,
$r(\sigma \{1,2\})$,..., $r(\sigma \{1,\dots , n\})$. For
$n\geq i\geq j$ let $r(\sigma :i:j)=r(\sigma \{1,\dots , i\}:j)$
in the notation of Section 1.3.

Note that $$\bar r(\sigma :i:j)=\bar r(\sigma (\{1,\dots , i\})
\bar r(\sigma (\{1,\dots , i-1\})\dots 
\bar r(\sigma (\{1,\dots , i-j+1\}).$$
Then $\{\bar r(\sigma :i_1:j_1)\dots \bar r(\sigma :i_k:j_k)
\ |\ ((i_1, j_1),\dots , (i_k, j_k))\in P^{\prime } \}$ 
spans  $\text {gr}\ Q_{n, \sigma }$
and so, by Theorem 1.3.8, is a basis for $\text {gr}\ Q_{n, \sigma}$.

Note that $\bar r(\sigma :i_1:j_1) \dots \bar r(\sigma :i_k:j_k)
\in G((i_1, j_1),\dots , (i_k, j_k))$. Denote this element by
$\bar r(\sigma :((i_1, j_1),\dots , (i_k, j_k)))$.

Now let 
$$x
=\sum _{((i_1, j_1),\dots , (i_s, j_s))\in P^{\prime }}
a_{((i_1, j_1),\dots , (i_s, j_s))}
\bar r(\sigma :((i_1, j_1),\dots , (i_s, j_s)))\in Q_{n, \sigma }
$$
and suppose $x\in \bigcap _{\tau \in S}Q_{n, \tau}$.
Then, by the $S_n$-invariance of the \break $G((i_1, j_1), \dots ,
(i_s, j_s))$, we have that 
$\bar r(\sigma :((i_1, j_1),\dots , (i_s, j_s)))$ is
$S$-invariant whenever $a_{((i_1, j_1),\dots , (i_s, j_s))}
\neq 0$.
This implies that $\bar r(\sigma :i_u:j_u)$ is $S$-invariant
for every $u$, $1\leq u\leq s$. So, $\{1,2,\dots , i_u\}$
is $S$-invariant, hence equal to $J_v$ for some $v$.
Then, writing $I_t$ for $I_t(S)$ and $J_t$ for $J_t(S)$, we have
 $$
\bar r(\sigma :i_u:j_u)=
\bar \Lambda _{j_u}(J_v)=
\bar \Lambda _{j_u}(I_1\cup\dots \cup I_v)=
\sum _{w_1+\dots +w_v=j_u} \bar \Lambda _{w_v}(I_v)
\dots \bar \Lambda _{w_1}(I_1).
$$
The result now follows.



One can say more about the subalgebra generated by
the coefficients of the polynomial $P(t)$.
Denote by $C$ the subalgebra of $Q_n$ generated by the 
elements $\Lambda _k=\tilde \Lambda _k(y_1,\dots, y_n)=
\Lambda _k(\{1,\dots, n\})$ for $k=1,\dots,n$.
\proclaim{Theorem 1.5.6} The algebra $C$ is
a free associative differential subalgebra of $Q_n$
generated by $\Lambda _k$, $k=1,2,\dots, n$. Also
$\theta (w)=w$ for and $w\in C$. 
\endproclaim
This follows from:
\proclaim{Proposition 1.5.7} For any $k=1,\dots,n$
$$ \partial \Lambda _k=(n-k+1)\Lambda _{k-1},$$
$$ \theta \Lambda _k=\Lambda _k.$$
\endproclaim
This shows that $C$ is a differential
subalgebra of the algebra $Q_n$. From \cite{GKLLRT} and
\cite{GR5} it follows that $C$ is also a free algebra. 

\medskip
              
\head 2. A map of $Q_n$ into a free skew-field
\endhead

\subhead 2.1  The free skew-field
$k<\!\!\!\!\!(x_1,\dots,x_n>\!\!\!\!)$ and a canonical derivation
\endsubhead

We will now recall the definition \cite{C1,C2} of
the free skew-field $k<\!\!\!\!\!(x_1,\dots,x_n>\!\!\!\!)$
generated by a set $\{x_1,\dots,x_n\}$.
Our purpose (Section 2.3) is to map $Q_n$ with
onto a subalgebra of $k<\!\!\!\!\!(x_1,\dots,x_n>\!\!\!\!)\ .$ 

The free associative algebra $k<x_1,\dots,x_n>$
on the set $\{x_1,\dots,x_n\}$ over a field $k$
has a universal field of fractions,
denoted $k<\!\!\!\!\!(x_1,\dots,x_n>\!\!\!\!)$ and called the 
{\it {free skew-field}} over $k$ on $\{x_1,\dots,x_n\}$.
The algebra $k<x_1,\dots,x_n>$ is naturally embedded into the algebra
$k<\!\!\!\!\!(x_1,\dots,x_n>\!\!\!\!)\ .$
The universality means that
if $D$ is any division ring and
$$\alpha : k<x_1,\dots,x_n> \longrightarrow D$$ is a 
homomorphism then there is a subring $R$ of 
$k<\!\!\!\!\!(x_1,\dots,x_n>\!\!\!\!)$ containing 
$k<x_1,\dots,x_n>$ and an extension of $\alpha $
to a homomorphism 
$$\beta: R \longrightarrow D$$ 
such that 
if $a \in R$ and $\beta (a) \neq 0$, then $a^{-1} \in R$.

Note that the symmetric group $S_n$ acts on
$k<x_1,\dots,x_n>$ by permuting subscripts and 
so $S_n$ acts on $k<\!\!\!\!\!\!(x_1,\dots,x_n>\!\!\!\!)\ .$
Denote the free skew-field \break $k<\!\!\!\!\!(x_1,\dots,x_n>\!\!\!\!)$ by $F$. 
The algebra $F$ has canonical partial derivations
$\partial _i, i=1,\dots,n$, such that

1)      $\partial _i(fg)=(\partial _if)g +
f(\partial _i)g$ for $f,g\in F$,

2)      $\partial _i(x_{j})=\delta _{ij}$ for $j=1,\dots,n$,

3)      $\partial _i(1)=0$.

Set $\nabla =\partial _1+...+\partial _n$. $\nabla $ is also 
a derivation, $\nabla (1)=0$ and $\nabla x_i=1$ for each $i$.
\smallskip
\subhead 2.2. Vandermonde quasideterminants \endsubhead 
Only now
are we going to prove that algebra $Q_n$ is correctly defined.
Instead of tedious work with defining relations (1.1) we will
use some quasideterminantal identities. These identities allow
us to construct a homomorphism of the algebra $Q_n$ onto a subalgebra in
$k<\!\!\!\!\!(x_1,\dots,x_n>\!\!\!\!)$. For this we will 
need a notion of the Vandermonde quasideterminant
developed in \cite {Gr3, GR5}. 
Let $R$ be a division algebra.
For $y_1,\dots, y_m \in R$ let 
$$\tilde V(y_1,\dots,y_m)$$
denote the Vandermonde matrix
$$
\left(\matrix y_1^{m-1}&...&y_m^{m-1}\\ y_1^{m-2}&...&y_m^{m-2}\\
...&...&...\\ y_1&...&y_m\\1&...&1\endmatrix\right).
$$
Then define 
$$V(y_1,\dots,y_m)=|\tilde V(y_1,\dots,y_m)|_{1m},$$ 
the corresponding Vandermonde quasideterminant \cite{GR3, GR5}.

Note that $V(y_1,\dots,y_m)$ is a rational function in $y_1,\dots,y_m$ 
and that it does not depend on the ordering of 
$y_1,\dots,y_{m-1}$ \cite{GR1-5}. 
Let elements $x_1,\dots,x_n$ belong to $R$ and assume that
the matrix $\tilde V(x_{i_1},\dots,x_{i_k})$ is invertible 
whenever $k=1,\dots,n$ and $i_1,\dots,i_k$ are distinct
integers, 
$1\leq i_1,\dots,i_k\leq n$. For any $i\neq i_1,\dots,i_k$ 
define the rational function
$x_{i_1...i_k, i}$ of $x_{i_1},\dots,x_{i_k}, x_i$
by the formula
$$ 
x_{i_1...i_k,i}=
V(x_{i_1},\dots,x_{i_k},x_i)
x_iV(x_{i_1},\dots,x_{i_k},x_i)^{-1}.
$$

The elements $x_{i_1...i_k,i}$ do not depend on the
ordering of $x_{i_1}, ..., x_{i_k}$.
\smallskip
\subhead 2.3 A subalgebra of a free
skew-field \endsubhead

Here we are going to construct a homomorphism of the 
quadratic algebra $Q_n$ into the free skew-field $F$.

For any $i$, set $x_{\emptyset, i}=x_i$. 
For any $A=\{i_1,\dots,i_k\}\subset \{1,\dots,n\}$,
$A\neq \emptyset$
and $i\notin A, 1\leq i \leq n$, set
$$
x_{A,i}=V(x_{i_1},\dots,x_{i_k},x_i)
x_iV(x_{i_1},\dots,x_{i_k},x_i)^{-1}. \tag 2.1
$$

Denote by $\hat Q_n$ the subalgebra of the free skew-field
$F$ generated by all elements $x_{A,i}$.    
\proclaim{Proposition 2.3.1} The subalgebra $\hat Q_n$ is 
a differential subalgebra of $(F,\nabla)$.   
\endproclaim
The proof follows from Lemma 2.3.2 and Corollary 2.3.3 below.
\proclaim{Lemma 2.3.2} For any subset of distinct numbers 
$\{i_1,\dots,i_s\}\subset \{1,\dots,n\}$
$$ \nabla V(x_{i_1},\dots,x_{i_s})=0. $$
\endproclaim
Lemma 2.3.2 and formula (2.1) imply the following:
\proclaim{Corollary 2.3.3} For all $A, i\notin A$
$$ \nabla x_{A,i}=1. $$
\endproclaim
 
\noindent{\bf Example}. $\nabla V(x_1,x_2)=\nabla (x_2-x_1)=0$
and so $\nabla x_{1,2}=1$.
\smallskip
Construct a homomorphism $\alpha $ of the algebra $Q_n$
into the skew-field $F$ by setting $\alpha (z_{A,i})=x_{A,i}$
for all $A, i\notin A$.
\proclaim{Conjecture 2.3.4} The homomorphism $\alpha $ is an embedding.
It defines an isomorphism of the differential quadratic
algebras $Q_n$ and $\hat Q_n$.
\endproclaim 
We can prove this conjecture for $n=2$.

\smallskip
\noindent {\bf Remark}. The derivation $\nabla $ transfers the identities 
(1.3b) into (1.3a) so the algebra $Q_n$ is defined by
the multiplicative identity (1.3b) and the derivation
$\partial $.             

The relations (1.1a) and (1.1b) imply 
more general rational relations between the generators
$z_{A,i}$ in the ``rational envelope" of the algebra $Q_n$. 
These relations follow from expressions for the $x_{A,i}$ 
described by Theorem 2.3.5 below. Let
$A\subset \{1,\dots,n\}$, 
$C\subset B\subset A$, $A\setminus B=\{i_1,\dots,i_p\}$ and 
$B\setminus C=\{j_1,\dots,j_q\}$.
\proclaim{Theorem 2.3.5} 
(a)\ Let $i\notin A$. Then
$$
x_{A,i}=V(x_{B,i_1},\dots,x_{B,i_p}, x_{B,i})
x_{B,i}V(x_{B,i_1},\dots,x_{B,i_p}, x_{B,i})^{-1}. \tag 2.2a
$$
(b)\ Let $j\notin B$. Then
$$
x_{C,j}=V(x_{B,j_1},\dots,x_{B,j_q}, x_{B,j})^{-1} 
x_{B,j}V(x_{B,j_1},\dots,x_{B,j_q}, x_{B,j}).\tag 2.2b
$$
\endproclaim

\noindent {\bf Examples}. a)\  If $A=\{i_1\}, B=\emptyset$ 
$$x_{i_1,i}=(x_{i_1}-x_i)x_i(x_{i_1}-x_i)^{-1}.$$

b)\  If $A=\{i_1,i_2\}, B=\{i_2\}$
$$x_{i_1i_2, i}=(x_{i_2,i}-
x_{i_2,i_1})
x_{i_2,i}(x_{i_2,i}-x_{i_2,i_1})^{-1}.$$

The last identity may be rewritten as
$$x_{i_2, i}=(x_{i_2,i}-x_{i_2,i_1})^{-1}
x_{i_1i_2,i}(x_{i_2,i}-x_{i_2,i_1}).$$

\noindent {\bf Remark}. Theorem 2.3.5 shows that the algebra 
$Q_n$ may be mapped not only into the free skew-field 
generated by the elements $x_k=x_{\emptyset ,k}$,
but into the free skew-field generated by the elements 
$x_{B,i}$ with $|B|=m$ for a given $m$.

\smallskip
\bigskip
\head 3. Noncommutative polynomials and their factorizations
\endhead
\smallskip
\subhead 3.1. Noncommutative polynomials \endsubhead 
Let $R$ be an associative algebra over a field $k$ of 
characteristic zero.  We denote 
by $x$ a noncommutative formal variable and by $t$ a 
commutative formal variable.  
We consider here the {\it associated} polynomials 
$\hat P(x)= a_0x^n + a_1x^{n-1} + ... +a_n$ and 
$P(t)=a_0t^n + a_1t^{n-1} + ... +a_n, \ a_0\neq 0$ over $R$.  
Recall that $x$ does not commute with the coefficients 
$a_0,\dots, a_n$ but $t$ is a commuting variable.
Relations between $\hat P(x)$ and $P(t)$ go back to Ore
\cite{O, L}; in particular one has the following:
\proclaim{Lemma 3.1.1} An element $\xi \in R$ is a root
of the polynomial $\hat P(x)$ if and only if 
$$ P(t)=Q(t)(t-\xi), $$
where $Q(t)$ is a polynomial over $R$. 
\endproclaim

A {\it generic} polynomial $\hat P(x)$ of degree $n$ over 
a division algebra has exactly $n$ roots.  
This follows from the next  result which is proved in \cite{BW}, 
see also \cite{L}.
\proclaim{Theorem 3.1.2} If a polynomial of degree $n$ over 
a division ring has  more than $n$ roots, then it has 
infinitely many roots.  
\endproclaim
Coefficients of polynomials with infinitely many roots were 
described in \cite{BW}, see also \cite{L}.

\noindent {\bf Remark}. For polynomials over a ring without division 
Theorem 3.1.2 is not true. A generic polynomial  of 
degree $m$ over the ring of  complex matrices of order 
$n$  has  $\binom {nm}n$ roots.

It was shown in \cite{GR5} that if a polynomial $\hat P(x)$
over a division algebra has $n$ roots in generic
position, then the associated polynomial $P(t)$ admits
a factorization into linear factors
$$P(t)=a_0(t-y_n)(t-y_{n-1})...(t-y_1). $$
Here we study all such decompositions and 
a quadratic algebra associated with them. 
\smallskip
\subhead 3.2. Factorizations of noncommutative polynomials
\endsubhead
A polynomial \break
$\hat P(x)= a_0x^n + a_1x^{n-1} + ... +a_n$ 
over a division algebra is called a 
{\it generic polynomial} if $\hat P(x)$ has exactly $n$ 
roots $x_1,\dots,x_n$ such that all rational expressions 
$x_{i_1...i_{k-1}, i_k}$ are defined and different from each other. 
A  polynomial $P(t)$  is called a generic polynomial 
if $\hat P(x)$ is a generic polynomial.
 
Set $x_{\emptyset, i}=x_i$. From \cite{GR3, GR5} we have:
\proclaim{Theorem  3.2.1} 
Let  $P(t)= a_0t^n + a_1t^{n-1} + ... +a_n$ be a generic polynomial.  
For any ordering $i_1,\dots,i_n$ of $\{1,\dots,n\},$
$$
P(t)=a_0(t-x_{i_1...i_{n-1},i_n})...(t-x_{i_1...i_{k-1}, i_k})...
(t-x_{\emptyset, i_1}). \tag 3.1
$$
Conversely, for any factorization 
$$
P(t)=a_0(t-\xi _n)(t-\xi _{n-1})...(t-\xi _1)
$$ 
there exists an ordering $i_1,i_2,\dots,i_n$ of  $\{1,2,\dots,n\},$ 
such that $\xi _k=x_{i_1...i_{k-1}, i_k}$ for $k=1,\dots,n$.
\endproclaim

{\it Proof}. In \cite{GR3, GR5} it was shown that
$$a_0^{-1}a_1=-(y_{i_1}+y_{i_2}+...+y_{i_n}), $$
$$a_0^{-1}a_2=\sum _{j<k}y_{i_k}y_{i_j}, $$
$$...$$
$$a_0^{-1}a_n=(-1)^ny_{i_n}y_{i_{n-1}}...y_{i_1}, $$ 
where $y_{i_k}=x_{i_1...i_{k-1}, i_k}$.

The factorization (3.1) is equivalent to this system. This proves the
first statement of the theorem.
The second statement follows from the first one and Theorem 3.1.2.

\proclaim{Corollary 3.2.2} There exist at most $n!$ factorizations 
of a generic polynomial of degree $n$ over a division algebra.    
\endproclaim

\noindent {\bf Example}. For a generic quadratic polynomial
$\hat P(x)=a_0x^2+a_1x+a_2$ and its roots $x_1, x_2$ one
has
$$x_{1,2}=(x_2-x_1)x_2(x_2-x_1)^{-1},$$
$$x_{2,1}=(x_1-x_2)x_1(x_1-x_2)^{-1},$$
$$P(t)=a_0(t-x_{1,2})(t-x_1)=a_0(t-x_{2,1})(t-x_2).$$
\smallskip
\subhead 3.3. Basic relations arising from factorizations of 
noncommutative polynomials
\endsubhead
 
Let $\hat P(x)= a_0x^n + a_1x^{n-1} + ... +a_n$ 
be a generic polynomial over a division algebra $R$ and 
let $x_1,\dots, x_n$ be its roots. We describe here basic relations 
for rational expressions in $x_{i_1i_2...i_k, i_{k+1}}$ for 
$k=0,\dots,n-1$.

Let $A=\{i_1,i_2,\dots,i_k\}$ be a subset of $\{1,2,\dots,n\}$ 
and $l\notin A$, $1\leq l \leq n$.  
\proclaim{Theorem 3.3.1} Let $|A|<n-1$. For any
$i,j\notin A$
$$x_{A\cup i, j}+x_{A,i}=x_{A\cup j, i}+x_{A,j}, \tag 3.2a $$
$$x_{A\cup i, j}x_{A,i}=
x_{A\cup j, i}x_{A,j}. \tag 3.2b $$
\endproclaim

\noindent {\bf Example}. Let $n=2$, $A=\emptyset$, $i=1, j=2$.
Then
$$x_{1,2}+x_1=x_{2,1}+x_2, $$
$$x_{1,2}x_1=x_{2,1}x_2.$$

Denote by $R_P$ the subalgebra of the algebra $R$ generated by all
$x_{A, i}$'s.
\proclaim {Corollary 3.3.2} There exists a unique epimorphism
$$\alpha :\ Q_n \to R_P $$
such that $\alpha (z_{A, i})=x_{A, i}$ for all $A$,
$i\notin A$.
\endproclaim

As in Theorem 2.3.5 more general relations for 
rational expressions in the $x_{A,i}$ are
described by the following theorem. Let
$A\subset \{1,\dots,n\}$, 
$C\subset B\subset A$, $A\setminus B=\{i_1,\dots,i_p\}$ and 
$B\setminus C=\{j_1,\dots,j_q\}$.
\proclaim{Theorem 3.3.3} 
(i)\ Let $i\notin A$. Then
$$
x_{A,i}=V(x_{B,i_1},\dots,x_{B,i_p}, x_{B,i})
x_{B,i}V(x_{B,i_1},\dots,x_{B,i_p}, x_{B,i})^{-1}. \tag 3.3a
$$
(ii)\ Let $j\notin B$. Then
$$
x_{C,j}=V(x_{B,j_1},\dots,x_{B,j_q}, x_{B,j})^{-1} 
x_{B,j}V(x_{B,j_1},\dots,x_{B,j_q}, x_{B,j}).\tag 3.3b
$$
\endproclaim

\noindent {\bf Examples}. a)\  If $A=\{i_1\}, B=\emptyset$ 
$$x_{i_1,i}=(x_{i_1}-x_i)x_i(x_{i_1}-x_i)^{-1}.$$

b)\  If $A=\{i_1,i_2\}, B=\{i_2\}$
$$x_{i_1i_2, i}=(x_{i_2,i}-
x_{i_2,i_1})
x_{i_2,i}(x_{i_2,i}-x_{i_2,i_1})^{-1}.$$

The last identity may be rewritten as
$$x_{i_2, i}=(x_{i_2,i}-x_{i_2,i_1})^{-1}
x_{i_1i_2,i}(x_{i_2,i}-x_{i_2,i_1}).$$
\bigskip
\head 4. Quadratic algebras associated with differential
polynomials \endhead
\subhead 4.1 Miura decompositions of differential
polynomials\endsubhead

In the next sections we transfer results obtained 
in previous sections to factorizations of differential polynomials.
Let $k$ be a field and $(R,D)$  a differential division algebra 
with unit. 
Here $D:R\rightarrow R$ is a $k$-linear map 
such that $D(ab)=(Da)b+a(Db)$ for $a,b\in R$. 
Sometimes we use the notation 
$Da=a^\prime $ and $D^ka=a^{(k)}$.
Consider an operator 
$L:R\rightarrow R$, $L=D^n+a_1D^{n-1}+...+a_n$, 
$a_i\in R$ for $i=1,\dots,n$. 
The action of $L$ on $R$ is given by the formula 
$L(\phi)=D^n\phi+a_1D^{n-1}\phi+...+a_n\phi$.
For $\phi _1,\dots,\phi _n\in R$ set
$$
\tilde W(\phi_1,\dots,\phi_m)=
\left(\matrix \phi_1^{(m-1)}&...&\phi_m^{(m-1)}\\
\phi_1^{(m-2)}&...&\phi_m^{(m-2)}\\ ..&...&..\\
\phi_1&...&\phi_m\endmatrix\right),
$$
and denote by $W(\phi_1,\dots,\phi_m)$ the quasideterminant 
$|\tilde W(\phi_1,\dots,\phi_m)|_{1m}$ of this matrix.
Suppose that $L$ has $n$ linearly independent solutions 
$\phi _i\in R$, i.e., $L(\phi _i)=0, i=1,\dots,n$. 
Suppose also that the matrix 
$\tilde W(\phi_{i_1},\dots,\phi_{i_m})$ is invertible for any set 
$\{i_1,\dots,i_m\}\subset \{1,\dots,n\}$.
Denote by $V$ the space of all solutions of $L$.
Let $F$ be the complete flag 
$\{F^{(1)}\subset F^{(2)}\subset ... \subset F^{(n)}\}$
such that each $F^{(k)}$ is generated by 
$\phi _{i_1},\dots,\phi _{i_k}$. Set
$$
b_k(F)=[DW(\phi _{i_1},\dots,\phi _{i_k})]
W(\phi _{i_1},\dots,\phi _{i_k})^{-1}.
$$

The elements $b_k(F)$ do not depend on the choice of the basis for $V$.
The following theorem was proved in \cite{EGR}.
\proclaim{Theorem 4.1.1}
$$
L= (D-b_n(F))(D-b_{n-1}(F))...(D-b_1(F)).\tag 4.1
$$
\endproclaim
It is easy to see that if $F_1,F_2$ are complete
flags in $V$ such that $F_1^{(m)}=F_2^{(m)}$
for $m=1,\dots,k$ and $F_1^{(k+2)}=F_2^{(k+2)},$
then:
\proclaim{Proposition 4.1.2}
$$
b_{k+2}(F_1)+b_{k+1}(F_1)=b_{k+2}(F_2)+b_{k+1}(F_2),
$$
$$
b_{k+2}(F_1)b_{k+1}(F_1)-Db_{k+1}(F_1)= 
b_{k+2}(F_1)b_{k+1}(F_2)-Db_{k+1}(F_2).
$$
\endproclaim
\smallskip
\subhead 4.2. Factorizations of differential
polynomials\endsubhead

Motivated by results from previous sections 
we are going to study differential algebras 
generated by the decompositions of differential polynomials. 
We do not assume here that these polynomials have any solutions.  

Let $(R,D)$ be a $k$-differential algebra with unit over 
a field $k$ and $L=D^n+a_1D^{n-1}+...+a_n$ 
a differential polynomial over $R$. 
For an element $g\in R$ set $u_p(g)=(D+g)^p(1)$.
For example, $u_0(g)=1$, $u_1(g)=g$, $u_2(g)=g^\prime +g^2$.
When $D=0$, $u_p(g)=g^p$.

Assume that the operator $L$ can be factorized as $L=L_i(D-f_i)$, 
where the $L_i$ are differential polynomials of degree 
$n-1$, $i=1,\dots,n$. 
Suppose that all square submatrices of the matrix 
$(u_p(f_q))$, $p=0,1,..,n-1$, $q=1,\dots,n$, are invertible. 
For distinct $i_1,\dots,i_m$, $m=1,\dots,n$, set 
$\theta (f_{i_1},\dots,f_{i_m})=|u_p(f_{i_s})|_{mm}$, 
where $p=0,\dots,m-1$, $s=1,\dots,m$. 
When $D=0$, the last quasideterminant is just a Vandermonde 
quasideterminant.

From a general property of quasideterminants \cite{GR1-GR5} 
it follows that $\theta (f_{i_1},\dots,f_{i_m})$ is 
symmetric in $i_1,\dots,i_{m-1}$.
Set
$$
f_{i_1,\dots,i_m}=\theta (f_{i_1},\dots,f_{i_m}) f_{i_m}
\theta (f_{i_1},\dots,f_{i_m})^{-1} + 
[D\theta (f_{i_1},\dots,f_{i_m})] \theta (f_{i_1},\dots,f_{i_m})^{-1}.
$$
The expressions $f_{i_1,\dots,i_m}$ are also symmetric in 
$i_1,\dots,i_{m-1}$.

\proclaim{Theorem 4.2.1} For any permutation 
$(i_1,\dots,i_n)$ of $\{1,\dots,n\}$
$$
L=(D-f_{i_1,...i_{n-1},i_n})...(D-f_{i_1,i_2})
(D-f_{i_1}). \tag 4.2
$$
\endproclaim

The decomposition (4.2) is similar to the decomposition (3.1)
of a polynomial $P(t)$.

\proclaim{Proposition 4.2.2} If $D\phi _i=f_i\phi _i$ for $i=1,\dots,n$, 
then decomposition (4.2) implies decomposition (4.1).
\endproclaim

\noindent {\bf Example}. For $n=2$
$$
L=(D-f_{1,2})(D-f_1)=(D-f_{2,1})(D-f_2),
$$
where 
$$
f_{1,2}=(f_2-f_1)f_2(f_2-f_1)^{-1}+
(f_2^\prime -f_1^\prime )(f_2-f_1)^{-1},
$$
$$
f_{2,1}=(f_1-f_2)f_1(f_1-f_2)^{-1}+
(f_1^\prime -f_2^\prime )(f_1-f_2)^{-1}.
$$

The following theorem generalizes formulas (3.2a,b) for 
factorizations of noncommutative polynomials. 

\proclaim{Theorem 4.2.3} Let $A\subset \{1,\dots,n\}$,
$|A|<n-1$. For any $i,j\notin A$
$$
f_{A\cup i,j}+f_{A,i}=f_{A\cup j,i}+f_{A,j}, \tag 4.3a
$$
$$
f_{A\cup i,j}f_{A,i}-f_{A,i}^\prime =
f_{A\cup j,i}f_{A,j}-f_{A,j}^\prime . \tag 4.3b
$$
\endproclaim
\noindent {\bf Example}. Let $n=2$, $A=\emptyset $, $i=1$,
$j=2$. Then
$$
f_{1,2}+f_1=f_{2,1}+f_2, 
$$
$$
f_{1,2}f_1-f_1^\prime =f_{2,1}f_2-f_2^\prime . 
$$

Formulas (4.3a,b) show that there exists a quadratic linear 
algebra $R_L$ associated with factorizations of the operator $L$. 
In fact, the algebra $R_L$ possesses a natural derivation $\nabla $.

\proclaim{Theorem 4.2.4} There exists a unique
derivation $\nabla :R_L\rightarrow R_L$
such that

i)      $\nabla f_{A,i}=1$ for any pair $A,i$,
$i\notin A$,

ii)     $\nabla D=D\nabla$.
\endproclaim

This theorem is a generalization of Corollary 2.3.3.
The relations between the $f_{A,i}$ for different pairs $A,i$ 
are given by the following statements.  
Let $A\subset \{1,\dots,n\}$, $B\subset A$,
$C\subset B$. Let 
$A\setminus B=\{i_1,\dots,i_p\}$,
$B\setminus C=\{j_1,\dots,j_q\}$.

\proclaim{Theorem 4.2.5} (i) If $i\notin A$, then
$$
f_{A,i}=\theta (f_{B,i_1},\dots,f_{B,i_p},f_{B,i})f_{B,i}
\theta (f_{B,i_1},\dots,f_{B,i_p},f_{B,i})^{-1}
$$
$$
+[D\theta (f_{B,i_1},\dots,f_{B,i_p},f_{B,i})] 
\theta (f_{B,i_1},\dots,f_{B,i_p},f_{B,i})^{-1}.
$$
(ii) If $j\notin B$, then   
$$
f_{C,j}=\theta (f_{B,j_1},\dots,f_{B,j_q},f_{B,j})^{-1}f_{B,j}
\theta  (f_{B,j_1},\dots,f_{B,j_q},f_{B,j})
$$
$$
-\theta (f_{B,j_1},\dots,f_{B,j_q},f_{B,j})^{-1}[D\theta (f_{B,j_1},\dots,f_{B,j_q},f_{B,j})].
$$
\endproclaim

This is a generalization of Theorem 3.3.3.

\enddocument